\documentclass[a4paper]{elsarticle}
\usepackage[utf8]{inputenc}
\usepackage{import}

\usepackage{fontawesome}

\usepackage{fullpage}

\usepackage{lipsum}

\usepackage{circledsteps}
\pgfkeys{/csteps/fill color=white}
\pgfkeys{/csteps/inner ysep=5pt}
\pgfkeys{/csteps/inner xsep=5pt}

\usepackage[english]{babel}

\usepackage{enumitem}
\setlist[1]{itemsep=-5pt}
\usepackage[margin=2.45cm]{geometry}
\usepackage{amsmath} 
\usepackage{amssymb,amsthm}
\usepackage{cases}
\usepackage{mathrsfs}
\usepackage{euscript}
\usepackage{bm}
\newtheorem{remark}{Remark}[section]
\usepackage{algorithm}
\usepackage{algpseudocode}
\theoremstyle{definition}

\usepackage{stmaryrd}
\usepackage{cases}

\usepackage{appendix}
\providecommand{\appendixname}{Appendix}
\usepackage{pdflscape}
\usepackage{tabularx}
\usepackage{booktabs}

\usepackage{multicol}
\usepackage{multirow}
\usepackage{xparse}
\usepackage{xspace}
\usepackage{makecell}
\usepackage{tablefootnote}

\usepackage{xcolor}
\usepackage{footnote} 
\usepackage{graphicx}
\usepackage{svg}
\usepackage{epstopdf}
\usepackage{tikz}
\usetikzlibrary{calc}
\usepackage{pgfplots}
\pgfplotsset{compat=newest}

\usepackage[percent]{overpic}
\usepackage{placeins} 

\usepackage[colorlinks=true]{hyperref}
\usepackage{cleveref}
\biboptions{sort&compress}
\usepackage[numbers]{natbib}
\usepackage{caption}

\usepackage[notref,notcite,final]{showkeys}

\usepackage[color=gray, backgroundcolor=yellow, textwidth=2cm, textsize=footnotesize]{todonotes}
\usepackage{setspace}
\usepackage{subcaption}
\usepackage{fullpage}
\usepackage{framed}
\usepackage[noprefix]{nomencl}

\makenomenclature

\renewcommand*\nompreamble{\begin{multicols}{2}}
\renewcommand*\nompostamble{\end{multicols}}
\usepackage{etoolbox}
\usepackage{mathtools}
\usepackage{empheq}
\usepackage{pdfpages}

\newcommand{\tr}{\mathrm{tr}\,}

\newcommand{\eps}{\bm{\varepsilon}}
\newcommand{\sig}{\bm{\sigma}}
\newcommand{\btau}{\bm{\tau}}
\newcommand{\bmu}{\bm{\mu}}
\newcommand{\gam}{\bm{\gamma}}
\newcommand{\zet}{\bm{\zeta}}

\newcommand{\setI}{\mathrm{I}}

\newcommand{\Ctens}{\bm{\mathsf{C}}}
\newcommand{\Dtens}{\bm{\mathsf{D}}}
\newcommand{\Etens}{\bm{\mathsf{E}}}
\newcommand{\Atens}{\bm{\mathsf{A}}}

\newcommand{\discrete}[1]{\mkern 1.2mu\underline{\mkern-1.2mu#1\mkern-1.2mu}\mkern 1.2mu}

\newcommand{\curr}{^{k+1}}
\newcommand{\past}{^{k}}

\newcommand{\setD}{\mathrm{D}}
\newcommand{\setE}{\mathrm{E}}
\newcommand{\setZ}{\mathrm{Z}}

\newcommand{\Be}{\mathbf{B}^{\varepsilon}_e}
\newcommand{\Bg}{\mathbf{B}^{\gamma}_e}
\newcommand{\Bz}{\mathbf{B}^{\zeta}_e}
\newcommand{\Nc}{\mathbf{N}^{\rchi}_e}

\newcommand{\Bet}{\mathbf{B}^{\varepsilon\,\mathrm{T}}_e}
\newcommand{\Bgt}{\mathbf{B}^{\gamma\,\mathrm{T}}_e}
\newcommand{\Bzt}{\mathbf{B}^{\zeta\,\mathrm{T}}_e}
\newcommand{\Nct}{\mathbf{N}^{\rchi\,\mathrm{T}}_e}

\DeclareMathOperator*{\arginf}{arg \, inf}
\DeclareMathOperator*{\sym}{sym}
\DeclareMathOperator*{\dive}{div}

\DeclareRobustCommand{\rchi}{{\mathpalette\irchi\relax}}
\newcommand{\irchi}[2]{\raisebox{\depth}{$#1\chi$}} 

\usepackage{stackengine}
\setstackgap{S}{1.75pt}

\newcommand{\trip}{\,\Shortstack{. . .}\,}

\makeatletter
\newcommand{\multiline}[1]{%
  \begin{tabularx}{\dimexpr\linewidth-\ALG@thistlm}[t]{@{}X@{}}
    #1
  \end{tabularx}
}
\makeatother

\definecolor{myblue}{RGB}{0, 0, 204}
\newcommand{\tR}[1]{{\textcolor{black}{#1}}}
\allowdisplaybreaks
\graphicspath 
{
  {./}
  {./source/}
}
\begin{document}

\begin{frontmatter}

\author[add1]{Jacinto Ulloa\corref{cor1}}
\ead{julloa@caltech.edu}
\author[add2]{Laurent Stainier}
\ead{laurent.stainier@ec-nantes.fr}
\author[add1]{Michael Ortiz}
\ead{ortiz@aero.caltech.edu}
\author[add1]{José E. Andrade\corref{cor1}}
\ead{jandrade@caltech.edu}

\tnotetext[t1]{{\itshape Postprint version.}}
\tnotetext[t2]{{\itshape Published version:}
       J.~Ulloa,  L.~Stainier, M.~Ortiz, and J.E.~Andrade.
       \newblock Data-driven micromorphic mechanics for materials with strain localization.
       \newblock {\em Computer Methods in Applied Mechanics and Engineering}, 2024. 
       {\itshape DOI: \tt\url{https://doi.org/10.1016/j.cma.2024.117180}} \vspace{1mm}}
       
\cortext[cor1]{Corresponding author}
\address[add1]{Division of Engineering and Applied Science, California Institute of Technology, Pasadena, CA 91125, USA}
\address[add2]{Nantes Universit{\'e}, Ecole Centrale Nantes, CNRS, GeM, UMR 6183, F-44000 Nantes, France}

\title{
Data-driven micromorphic mechanics for materials with strain localization\tnoteref{t1,t2}
}

\begin{abstract}
This paper explores the role of generalized continuum mechanics, and the feasibility of model-free data-driven computing approaches thereof, in solids undergoing failure by strain localization. Specifically, we set forth a methodology for capturing material instabilities using data-driven mechanics without prior information regarding the failure~mode. We show numerically that, in problems involving strain localization, the standard data-driven framework for Cauchy/Boltzmann continua fails to capture the length scale of the material, as expected. We address this shortcoming by formulating a generalized data-driven framework for micromorphic continua that effectively captures both stiffness and length-scale information, as encoded in the material data, in a model-free manner. These properties are exhibited systematically in a one-dimensional softening bar problem and further verified through selected~plane-strain~problems.
\end{abstract}


\end{frontmatter}

\section{Introduction}

It is well-acknowledged that computational models based on classical Cauchy/Boltzmann continua strongly deviate from reality when microstructural effects modulate the material response. This is evident in the case of solids exhibiting different forms of strain localization, a ubiquitous process resulting from material instability, often related to softening responses or non-associative evolution laws. Examples include shear bands in granular materials driven by complex particle interactions~\cite{alonso2005,ando2012,karapiperis2021}, damage and fracture in quasi-brittle materials governed by evolving microcracks~\cite{bazant1986,labuz2001,zhao2018}, and ductile failure in metals emerging from the interplay between dislocations and void growth~\cite{asaro1977,cuitino1996,tvergaard2004}. \tR{The inadequacy of classical theories to describe such behaviors arises from their simple nature, preserving a continuum representation of the material at arbitrarily small scales. The de facto approach to this problem in the mechanics community is to introduce, in one way or another, length scales that reflect the heterogeneities of the microstructure.} Aside from strain localization, this approach also allows the description of size effects~\cite{bazant1999,ortiz2000,gao2001,forest2000} and wave dispersion phenomena~\cite{madeo2015,misra2016,dayal2017}.

From a computational perspective, the central dilemma of standard local models with strain localization is the loss of ellipticity of the boundary value problem (BVP), leading to pathological mesh dependence in finite element simulations with vanishing energy dissipation. This issue emerged in the 1970s and has since been addressed by constitutive models on various fronts. A conventional approach involves endowing internal variables with non-local effects. A substantial body of work exists in this context for gradient plasticity~\cite{muhlhaus1991,hutchinson1997,gurtin2009,bardella2010,martinez2016,scherer2021,ariza2023}, gradient damage~\cite{peerlings2004,comi1999,marigo2016}, and integral-type non-local constitutive models~\cite{pijaudier1987,bazant2002,grassl2006plastic}. Furthermore, non-local effects emerge naturally in fracture modeling via phase-field approaches~\cite{bourdin2000,marigo2016}, which can be equipped with additional regularization mechanisms~\cite{miehe2016,miehe2017,rodriguez2018,lancioni2020,suh2020,kristensen2020,ulloa2021b}, and peridynamics~\cite{silling2000,lipton2014,bhattacharya2023}.

It is noteworthy that before these developments, the need for generalized continuum theories for microstructured media had long been recognized, starting with the seminal framework of \emph{micropolar} continua developed by the Cosserats~\cite{cosserat1909}. Therein, a continuum material point is endowed with independent rotations, associated with the microstructure, whose gradient enters the formulation and serves as a regularization measure. The more general theory of \emph{micromorphic} continua was later elaborated in the works of~\citet{eringen1964}, \citet{mindlin1964}, and~\citet{germain1973}, incorporating a full, generally non–symmetric, micro-deformation tensor as an additional set of degrees of freedom. This framework recovers several models as special cases; for instance, the micropolar theory is recovered when the micro-deformation tensor is purely skew-symmetric, while the second-gradient theory is recovered when the micro-deformation equals the compatible displacement gradient. Moreover, the notion of micromorphic degrees of freedom has been extended to plasticity~\cite{forest2003,regueiro2009,bryant2019,rys2020,lindroos2022} and damage~\cite{aslan2011,miehe2017,brepols2017,yin2022} variables. For a thorough treatise of micromorphic continuum theories, we refer the reader to~\citet{forest2009} and~\citet{neff2014}.

All the references above are examples of continuum theories with constitutive models. In particular, the classical stress-strain relation in Cauchy continua--a phenomenological equation itself--is complemented by generalized stress-strain relations. Depending on the modeling assumptions, the number of additional elasticity parameters may increase significantly. Moreover, modeling inelastic behaviors, such as plasticity or damage, further requires an additional set of parameters and, perhaps more importantly, constitutive functions and evolution equations tailored to match experimental results. The inevitable caveat of this conventional modeling paradigm is epistemic uncertainty, particularly for complex material behavior.

\tR{The work of~\citet{kirchdoerfer2016} on data-driven computing presents a novel perspective that aims to resolve this shortcoming of constitutive modeling in a general context. This framework enforces the fundamental balance laws of continuum mechanics, as in classical problems, but extracts the material behavior directly from empirical data rather than from a constitutive model. The method works for both elastic~\cite{kirchdoerfer2016,conti2018data} and inelastic materials~\cite{eggersmann2019,ciftci2022model} and has been extended to finite deformations~\cite{nguyen2018data,conti2020data,kuang2023data}, composite structures~\cite{xu2020data}}, dynamics~\cite{kirchdoerfer2018data,garcia2023}, viscoelasticity~\cite{salahshoor2022model}, fracture mechanics~\cite{carrara2020data}, and molecular dynamics~\cite{bulin2023}. Data-driven computing has also been used in a multi-scale setting to capture the response of real granular materials~\cite{karapiperis2021a}, including breakage mechanics in crushable sand~\cite{ulloa2023}. All these works are based on Cauchy continua and are thus unable to handle material instabilities. An exception is found in~\citet{karapiperis2021b}, where data-driven computing is extended to micropolar media, a framework suitable for shear-banding driven by significant micro-rotations. \tR{Moreover, an idea for data-driven analysis using strain gradients is given in~\citet{kamasamudram2023strain} based on Lipschitz regularization.}

The present study presents a generalized data-driven mechanics framework for micromorphic continua. The motivation for this framework stems from softening materials with strain localization, where standard data-driven computing yields mechanically admissible states that fail to capture the length scale of the material, i.e., the expected spatial distribution of the state variables. This limitation and its resolution via micromorphic mechanics are showcased using synthetic data in a 1D softening bar governed by a gradient damage model. The framework is further validated in 2D problems considering damage and non-associative plasticity. Remarkably, no a priori data assignment is required for the proposed framework to capture the reference response in the presence of strain localization.

The paper is structured as follows. Section~\ref{sec:micromorphic} gives an overview of micromorphic continuum mechanics for the reader's convenience. Then, section~\ref{sec:DD} presents the proposed data-driven formulation and numerical implementation. Section~\ref{sec:num} shows the numerical examples and section~\ref{sec:conclusion} provides concluding remarks. \tR{Finally,~\ref{sec:app} presents the reference constitutive models used to generate the synthetic data sets, while~\ref{sec:app_data} outlines possibilities for data acquisition in practical scenarios.}

\section{Overview of generalized micromorphic continua} 
\label{sec:micromorphic} 

\subsection{Problem setting}

\noindent Consider a solid body $\Omega\subset\mathbb{R}^n$, $n\in\{1,2,3\}$, evolving in a time interval $\setI=[0,T]$. In a classical Cauchy theory with linearized kinematics, the motion of the continuum solid is fully characterized by the macro-displacement field $\bm{u}:\Omega\times\setI\to\mathbb{R}^n$. Conversely, \tR{the theory of micromorphic continua~\cite{eringen1964,mindlin1964}} further considers that attached to each macro-coordinate $\bm{x}\in\Omega$ is a \emph{micro-continuum} body $\Omega^\mathrm{m}$, subject to independent deformations and defined through micro-coordinate axes $\xi_i$ whose origin moves with the macro-displacement (figure 1). Then, the motion of the microstructured solid is given by a \emph{micro-displacement} field, conveniently expressed as the sum of the macro-displacement and a fluctuation~term: 
$$\bm{v}(\bm{x},\bm{\xi})=\bm{u}(\bm{x})+\bm{\rchi}(\bm{x})\cdot\bm{\xi}\,, \quad \tR{\rchi_{ij}=\partial v_{i}/\partial\xi_j}\,.$$
Evidently, in this form, the micro-deformation field $\bm{\rchi}:\Omega\times\setI\to\mathbb{R}^{n\times n}$ is homogeneous in $\Omega^\mathrm{m}$ and thus independent of~$\bm{\xi}\in\Omega^\mathrm{m}$, resulting in a micromorphic theory of degree 1~\cite{germain1973}. 

\begin{figure}[h!]
    \vspace{1em}
    \centering
    \small
    \vspace{1em}
    \includeinkscape[scale=0.8]{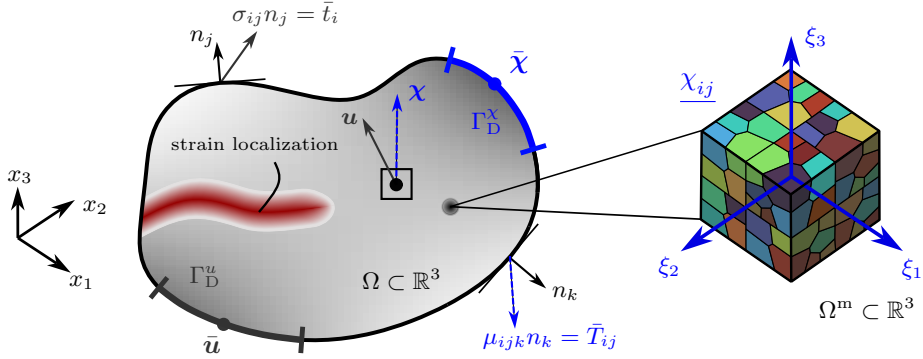}
    \vspace{1em}
    \caption{Schematic representation of a BVP in a micromorphic continuum (left) with a generic microstructure (right). The dashed blue arrows represent second-order tensors. The first-order theory adopted in this work assumes that the micro-deformation ${\rchi}_{ij}$ is homogeneous in $\Omega^\mathrm{m}$ but varies in $\Omega$.}
    \label{fig:micbvp}
\end{figure}

Consequently, we may formulate a BVP at the macro-scale in terms of the primary fields $\bm{u}$ and $\bm{\rchi}$, with the latter representing an additional set of degrees of freedom (DOFs). The solid boundary $\Gamma$ is then decomposed into a Dirichlet part $\Gamma^u_\mathrm{D}$ with imposed displacements $\bar{\bm{u}}:\Gamma^u_\mathrm{D}\times\setI\to\mathbb{R}^n$ and a Neumann part $\Gamma^u_\mathrm{N}$ with applied forces $\bar{\bm{t}}:\Gamma^u_\mathrm{N}\times\setI\to\mathbb{R}^n$, with $\Gamma^u_\mathrm{D}\cup\Gamma^u_\mathrm{N}=\Gamma$ and $\Gamma^u_\mathrm{D}\cap\Gamma^u_\mathrm{N}=\emptyset$. Similarly, the boundary conditions for $\bm{\rchi}$ consist of imposed micro-deformations $\bar{\bm{\rchi}}:\Gamma^\rchi_\mathrm{D}\times\setI\to\mathbb{R}^{n\times n}$ and applied double forces $\bar{\bm{T}}:\Gamma^\rchi_\mathrm{N}\times\setI\to\mathbb{R}^{n\times n}$, with $\Gamma^\rchi_\mathrm{D}\cup\Gamma^\rchi_\mathrm{N}=\Gamma$ and $\Gamma^\rchi_\mathrm{D}\cap\Gamma^\rchi_\mathrm{N}=\emptyset$.\footnote{We note that the physical significance of Dirichlet boundary conditions on $\bm{\rchi}$ is debatable.} For the sake of conciseness, volume forces and volume double forces are omitted, without losing generality.

\subsection{Kinematics and balance laws}

The micromorphic continuum of degree 1 involves the following generalized strain measures:
\begin{align}
&\text{Compatible strain} &\eps=\sym\,\nabla\bm{u}\,,&&& &\varepsilon_{ij}=\frac{1}{2}(u_{i,j}+u_{j,i})\,;&&&
\label{eq:epsilon}   \\
&\text{Relative strain}   &\gam=\nabla\bm{u} - \bm{\rchi}\,,&&& &\gamma_{ij}=u_{i,j}-\rchi_{ij}\,;&&&
\label{eq:gamma} \\ 
&\text{Micro-deformation gradient} &\zet=\nabla\bm{\rchi}\,,&&& &\zeta_{ijk}=\rchi_{ij,k}\,.&&&
\label{eq:zeta}
\end{align}
The power density of internal mechanical forces reads
\begin{equation}
p^\mathrm{int}=\sig:\dot{\eps}+\btau:\dot{\gam}+\bmu\trip\dot{\zet}\,, 
\label{eq:pdens}
\end{equation}
where $\sig$ is the second-order, symmetric Cauchy stress tensor; $\btau$ is the second-order, generally non-symmetric relative stress tensor; and $\bmu$ is the third-order double-stress tensor. The principle of virtual power
\begin{equation}
\int_\Omega\big(\sig:\dot{\eps}+\btau:\dot{\gam}+\bmu\trip\dot{\zet}\big)\,\mathrm{d}\bm{x}-\int_\Gamma\big(\bar{\bm{t}}\cdot\dot{\bm{u}} + \bar{\bm{T}}:\dot{\bm{\rchi}}\big) \,\mathrm{d}S=0
\label{eq:pvp}
\end{equation}
then yields the following balance equations $\forall\,t\in\setI$:
\begin{align}
&\text{Stress equilibrium} \qquad &\dive\,(\sig+\btau)=\bm{0} \ \ \ \text{in} \ \ \ \Omega\,,&&& \quad &(\sig+\btau)\cdot\bm{n}=\bar{\bm{t}}  \ \ \ \text{on} \ \ \ \Gamma_{\mathrm{N}}^{u}\,;&&&
\label{eq:eqsig}   \\
&\text{Double-stress equilibrium} \qquad  &\dive\,\bmu+\btau=\bm{0} \ \ \ \text{in} \ \ \ \Omega\,,&&& \quad &\bmu\cdot \bm{n}=\bar{\bm{T}}  \ \ \ \text{on} \ \ \ \Gamma_{\mathrm{N}}^{\rchi}\,.&&&
\label{eq:eqmu}
\end{align}
This system must be solved for the primary fields $\bm{u}\in\mathrm{H}^1(\Omega;\mathbb{R}^n)$ and $\bm{\rchi}\in\mathrm{H}^1(\Omega;\mathbb{R}^{n\times n})$ subject to the boundary conditions $\bm{u}=\bar{\bm{u}}$ on $\Gamma^u_\mathrm{D}$ and $\bm{\rchi}=\bar{\bm{\rchi}}$ on $\Gamma^\rchi_\mathrm{D}$.\footnote{When $\gam\equiv\bm{0}$, we recover the second-gradient theory with a single balance equation in terms of $\bm{u}\in\mathrm{H}^2(\Omega;\mathbb{R}^n)$.} Clearly, at this point, we require relations between the generalized stresses $\{\sig,\btau,\bmu\}$ and the conjugate generalized strains $\{\eps,\gam,\zet\}$. 

\subsection{Material behavior}
\label{sec:mmorph_material}

\subsubsection{Elasticity}

\noindent  As in classical Cauchy continua, conventional theories for generalized continua invoke constitutive models to establish the relation between $\{\sig,\btau,\bmu\}$ and $\{\eps,\gam,\zet\}$. In the case of linear elasticity, a linear mapping $(\eps,\gam,\zet)\mapsto \{\sig,\btau,\bmu\}(\eps,\gam,\zet)$ is employed, usually in terms of a \tR{quadratic Helmholtz energy density $\psi\coloneqq{\psi}(\eps,\gam,\zet)$, such that}
\begin{equation}
\sig=\frac{\partial\psi}{\partial\eps}\,, \qquad \btau=\frac{\partial\psi}{\partial\gam}\,, \qquad  \bmu=\frac{\partial\psi}{\partial\zet}\,.\label{eq:consmm}
\end{equation}
A fairly general form can be written as
\begin{equation}
\tR{{\psi}(\eps,\gam,\zet)\coloneqq\frac{1}{2}\big(\eps:\Ctens:\eps+\gam:\Dtens:\gam+\eps:\Etens:\gam+\zet\trip\Atens\trip\zet\big)\,,}
\label{eq:fremm}
\end{equation}
where $\Ctens$, $\Dtens$, $\Etens$ (4th-order), and $\Atens$ (6th-order) are generalized elasticity tensors. For isotropic materials, these tensors involve up to 18 parameters~\cite{mindlin1964}, most of which are difficult to obtain experimentally. Of course, this general model includes several special cases with fewer parameters (table~\ref{tab:mmorph}). \tR{For instance, the micropolar model, with $\bm\rchi$ purely skew-symmetric, preserves 6~parameters. \tR{Moreover, the coupling tensor $\Etens$ is often assumed to vanish for the sake of simplicity~\cite{neff2014,broese2018}}, leaving 16 parameters in the general theory. Then, the microstrain model, with $\bm\rchi$ purely symmetric, is reduced from 11 to 9 parameters. On the other hand, employing the curl of the microdeformation rather than the full gradient, the relaxed micromorphic model of~\citet{neff2014} is reduced from 9 to 7 parameters.} 

\begin{remark}
In previous works~\cite{forest2004, dillard2006,broese2018}, the following forms are adopted for the microstrain simplification of the general theory: $\mathsf{C}_{ijkl}=\lambda\delta_{ij}\delta_{kl}+G(\delta_{ik}\delta_{jl}+\delta_{il}\delta_{jk})$, where $\lambda$ and $G$ are the Lamé coefficients; $\mathsf{D}_{ijkl}=c_1\mathsf{C}_{ijkl}$, with $c_1\geq 1$; and $\mathsf{A}_{ijklmn}=c_1\ell_\chi^2\mathsf{C}_{ijlm}\delta_{kn}$, where $\ell_\chi\geq0$ is a length scale parameter. A micromorphic continuum of this type, with only 4 elasticity parameters, is considered for the sake of simplicity in the models presented in~\ref{sec:app} and employed in the numerical simulations of section~\ref{sec:num}. 
\label{rem_1}
\end{remark}

\begingroup
\begin{table}[]
\centering
\small
\caption{Some special cases of the first-order micromorphic theory for isotropic elastic materials with a free energy of the form~\eqref{eq:fremm}. See~\citet{forest2006} and~\citet{neff2014} for thorough discussions.}
\renewcommand{\arraystretch}{1.25} 
\begin{tabular}{lccc}
Model & Special case & No. DOFs & No. Params. \\ \hline
Cauchy                 & $\bm{\rchi}\rightsquigarrow\bm{0}$      &  3 $(\bm{u})$ & 2 \\
Second-gradient~\cite{germain1973b}        & $\bm{\rchi}\rightsquigarrow\bm{\eps}$   &  3 $(\bm{u})$ & 3 \\
Micropolar~\cite{cosserat1909}             & $\bm{\rchi}\rightsquigarrow\mathrm{skew}\,\bm{\rchi}$ & 6 $(\bm{u},\mathrm{skew}\,\bm{\rchi})$ & 6 \\
Microstrain~\cite{forest2006} & $\bm{\rchi}\rightsquigarrow\mathrm{sym}\,\bm{\rchi}$  & 9 $(\bm{u},\sym\bm{\rchi})$ & 11 \\
\multirow{2}{*}{Relaxed microm.~\cite{neff2014}} &   $\hspace*{-0.28cm} \nabla\bm{\rchi}\rightsquigarrow\mathrm{curl}\,\bm{\rchi}$  & \multirow{2}{*}{12 $(\bm{u},\bm{\rchi})$} & \multirow{2}{*}{9}  \\
&    $ \gam\rightsquigarrow\sym\,\bm{\gam}$  & 
\end{tabular}
\label{tab:mmorph}
\end{table}
\endgroup

\subsubsection{Inelasticity}

Extensions of generalized continuum theories to plasticity, damage, and other inelastic phenomena~\cite{forest2009} typically resort to the framework of thermomechanics with internal variables. In this setting, we introduce a generic set of internal variables $\mathbf{q}$ and consider free energies depending on both $\mathbf{q}$ and (possibly) its gradient, i.e., $\psi\coloneqq{\psi}(\eps,\gam,\zet,\mathbf{q},\nabla\mathbf{q})$. The internal variables obey evolution laws,  generally stated as differential inclusions of the form
\begin{equation}
\frac{\delta\psi}{\delta\mathbf{q}}(\eps,\gam,\zet,\mathbf{q},\nabla\mathbf{q})\; +\; \partial{\phi}(\dot{\mathbf{q}})\ni\bm{0}\,.
\label{eq:evol}
\end{equation}
Here, $\phi\coloneqq{\phi}(\dot{\mathbf{q}})\geq 0$ is a convex dissipation potential complying with the second law of thermodynamics. Evidently, the constitutive relations~\eqref{eq:consmm} are now functions of $\mathbf{q}$ and thus become history-dependent. Models of this type can be derived rigorously from the principle of virtual power~\cite{maugin1980} or, more generally, from the energetic formulation~\cite{ortiz1999,mielke2015}, as shown for a gradient damage model in~\ref{sec:app1}. It is worth noting that the dissipation potential may also depend on the current state, which allows us to consider non-associative models~\cite{francfort2018,ulloa2021}. \ref{sec:app2} presents a plasticity model of this type.

\subsection{Data-driven approach}
\label{sec:mmorph_dd}

In the standard formulation outlined above, even classical Cauchy continuum models are susceptible to epistemic uncertainty, particularly when dealing with complex inelastic behavior. This challenge is further amplified in the framework of generalized continua, requiring additional constitutive functions and parameters. Thus, in this study, we shift our focus towards the realm of data-driven mechanics~\cite{kirchdoerfer2016} for inelastic materials~\cite{eggersmann2019} and extend this approach to encompass generalized micromorphic continua. In this setting, instead of relying on constitutive relations and evolution equations such as~\eqref{eq:consmm} and~\eqref{eq:evol}, we assume that the material response in terms of generalized strains $\{\eps,\gam,\zet\}$ and stresses $\{\sig,\btau,\bmu\}$ is encapsulated in a finite data set, presumably obtained from experimental measurements or high-fidelity numerical simulations. We shall then make direct use of this data to solve the fundamental equations~\eqref{eq:epsilon}--\eqref{eq:zeta} and~\eqref{eq:eqsig}--\eqref{eq:eqmu}, without resorting to an empirical model for closure.

Of course, embarking on this path presupposes the availability of data sets comprising generalized coordinates $\{(\eps,\gam,\zet),(\sig,\btau,\bmu)\}$. Measurements of compatible strains and the corresponding Cauchy stresses, $(\eps,\sig)$, may be obtained, for instance, from high-fidelity virtual simulations and standard homogenization theory, as done in recent works on granular media~\cite{karapiperis2021a,gorgogianni2023,ulloa2023}. Acquiring data for the generalized quantities $(\gam,\btau)$ and $(\zet,\bmu)$ poses an additional challenge. In this context, a data-driven framework for micropolar continua using micro-rotation and couple-stress data directly extracted from LS-DEM simulations of sand has been recently presented~\cite{karapiperis2021b}. \tR{Extending virtual measurements to encompass generalized micromorphic continua appears feasible in light of previous works on this topic (see~\ref{sec:app_data}).}

In this study, we focus not on data acquisition but on developing a data-driven framework that incorporates generalized stress-strain data to simulate materials with strain localization. Consequently, we limit our studies in the sequel to synthetic data obtained from reference micromorphic constitutive models (\ref{sec:app}), \tR{while possibilities for identifying material data from other sources in practical scenarios are outlined in \ref{sec:app_data}}. The following section presents the proposed data-driven formulation.

\section{Set-oriented data-driven formulation of micromorphic continua}
\label{sec:DD}

This section presents the extension of data-driven computing to micromorphic materials. In agreement with most references on data-driven mechanics, we express the formulation in a discretized form, i.e., considering a finite-dimensional system. Consequently, the fundamental kinematic relations~\eqref{eq:epsilon}--\eqref{eq:zeta} and balance laws~\eqref{eq:eqsig}--\eqref{eq:eqmu} are restated in a finite element setting. The resulting equations are employed in the statement of the data-driven problem that follows.

\subsection{Discrete kinematics and balance equations}
\label{sec:fem}

Let us consider hereafter an isotropic solid discretized into $N$ nodes and $M$ material points, leading to the assembly of nodal displacement and external force vectors denoted as $\discrete{\mathbf{u}}=\{{\bm{u}}_i\in\mathbb{R}^{n}\}_{i=1}^N$ and $\discrete{\mathbf{f}}^\mathrm{ext}=\{{\bm{f}}^\mathrm{ext}_i\in\mathbb{R}^{n}\}_{i=1}^N$, respectively. Similarly, the nodal micro-deformations and double forces are collected in vectors denoted as $\discrete{\bm{\rchi}}=\{{\bm{\rchi}}_i\in\mathbb{R}^{n_\chi}\}_{i=1}^N$ and $\discrete{\mathbf{m}}^\mathrm{ext}=\{{\bm{m}}^\mathrm{ext}_i\in\mathbb{R}^{n_\chi}\}_{i=1}^N$, respectively. Here, the size of the nodal micro-deformation vector, $n_\chi$, depends on the symmetry conditions. For instance, $n_\chi=n^2$ for the full first-order micromorphic model and $n_\chi=n(n+1)/2$ for the microstrain simplification (see table~\ref{tab:mmorph}). 

On the other hand, the time interval $\setI$ is discretized into a finite number of time steps $0 =t^0<\dots < t^k < t^{k+1}< \dots < t^K = T$. Hereinafter, unless required for clarity, we denote a quantity at a given time $t\past$ as $\Box\past\coloneqq\Box(t\past)$ and omit the superscript at the next time step, so that $\Box\coloneqq\Box\curr$.

In this section, second- and higher-order tensors evaluated at a point in space are expressed in their Voigt representation.  Then, the discrete form of the kinematic relations~\eqref{eq:epsilon}--\eqref{eq:zeta} reads
\begin{align}
&\eps_e = \Be\,\discrete{\mathbf{u}}\;+\;\bm{g}^u_e\,,
\label{eq:epsilon_disc}   \\
&\gam_e = \Bg\,\discrete{\mathbf{u}}\;-\;\Nc\,\discrete{\bm{\rchi}}\,,
\label{eq:gamma_disc} \\ 
&\zet_e = \Bz\,\discrete{\bm{\rchi}}\;+\;\bm{g}^\rchi_e\,,
\label{eq:zeta_disc}
\end{align}
where $\bm{g}^u_e\in\mathbb{R}^{n_\varepsilon}$ and $\bm{g}^\rchi_e\in\mathbb{R}^{n_\zeta}$ follow from Dirichlet boundary conditions, with $n_\varepsilon=n(n+1)/{2}$ and $n_\zeta=n_\chi(n_\chi+1)/2$. Similarly, the balance equations~\eqref{eq:eqsig}--\eqref{eq:eqmu} take the discrete forms
\begin{align}
&\sum_{e=1}^{M}w_e\Big[\Bet\,\sig_e\;+\;\Bgt\,\btau_e\Big]-\discrete{\mathbf{f}}^\mathrm{ext}=\bm{0}\,,
\label{eq:eqsig_disc}\\
&\sum_{e=1}^{M}w_e\Big[\Bzt\,\bmu_e\;-\;\Nct\,\btau_e\Big]-\discrete{\mathbf{m}}^\mathrm{ext}=\bm{0}\,,
\label{eq:eqmu_disc}
\end{align}
where $w_e$ is an integration weight. In these expressions, $\mathbf{N}^\rchi_e\in\mathbb{R}^{n_\chi\times N n_\chi}$ is the shape function matrix for the nodal micro-deformation vector $\bm{\discrete{\rchi}}$ while $\Be\in\mathbb{R}^{n_\varepsilon \times N n }$, $\Bg\in\mathbb{R}^{n_\chi\times N n_\chi}$, and $\Bz\in\mathbb{R}^{n_\zeta\times N n_\chi}$ are standard discrete differential operators.

\subsection{The data-driven problem}

The data-driven problem requires a set-oriented reformulation of the discrete system described above. In this setting, the state of the solid at a material point $e$ is defined in terms of generalized coordinates $\mathbf{z}_e$ in a local phase space $\setZ_e$. For the present micromorphic formulation, we have, in general,
$$\mathbf{z}_e\coloneqq(\eps_e,\gam_e,\zet_e,\sig_e,\btau_e,\bmu_e)\in\mathrm{Z}_e,\quad \mathrm{Z}_e\coloneqq\mathbb{R}^{n_\varepsilon}\times\mathbb{R}^{n_\chi}\times\mathbb{R}^{n_\zeta}\times\mathbb{R}^{n_\varepsilon}\times\mathbb{R}^{n_\chi}\times\mathbb{R}^{n_\zeta}\,.$$
The local phase space is endowed with the following metric (cf.~\eqref{eq:fremm} with $\Etens=\bm{0}$):
\begin{equation}
\Vert \mathbf{z}_e \Vert^2 \coloneqq \frac{1}{2}\Big(\Ctens_e\eps_e\cdot\eps_e + \Ctens_e^{-1}\sig_e\cdot\sig_e + \Dtens_e\gam_e\cdot\gam_e + \Dtens_e^{-1}\btau_e\cdot\btau_e + \Atens_e\zet_e\cdot\zet_e + \Atens_e^{-1}\bmu_e\cdot\bmu_e\Big)\,.
\label{eq:locmet}
\end{equation}
Here, $\Ctens_e\in\mathbb{R}^{n_\varepsilon\times n_\varepsilon}_{\mathrm{sym},+}$, $\Dtens_e\in\mathbb{R}^{n_\chi\times n_\chi}_{\mathrm{sym},+}$, and $\Atens_e\in\mathbb{R}^{n_\zeta\times n_\zeta}_{\mathrm{sym},+}$ are positive-definite matrices of a  purely numerical nature, in principle, unrelated to the material behavior. Globally, the local coordinates $\mathbf{z}\coloneqq\{\mathbf{z}_e\in\mathrm{Z}_e\}_{e=1}^M\in\mathrm{Z}$ represent points in the global phase space $\mathrm{Z}=\mathrm{Z}_1\times\dots\times\mathrm{Z}_M$. In view of the discrete formulation of section~\ref{sec:fem}, the global distance metric reads
	\begin{equation}
	    \Vert \mathbf{z} \Vert^2 \coloneqq \sum_{e=1}^M w_e \Vert \mathbf{z}_e \Vert^2\,.
	    \label{eq:globmet}
	\end{equation}
\tR{Setting the metric tensors at the same order of magnitude of the elasticity parameters of a reference material has shown to yield good convergence behavior~\cite{eggersmann2021}. However, using the secant moduli, identifiable from the data, seems to give better results during softening stages  (section~\ref{sec:num}). We further note that a more objective choice of metric tensors consists of evaluating optimal moduli~\cite{karapiperis2021a}.}

We proceed to define a subset of $\mathrm{Z}$ containing all mechanically admissible states:
	\begin{equation}
	\setE\curr\coloneqq\Big\{\mathbf{z}\in\mathrm{Z} \ \colon \  \text{\eqref{eq:epsilon_disc}--\eqref{eq:eqmu_disc} } \ \text{at } \ t=t\curr\Big\}\,,
	\label{eq:Eset}
	\end{equation}
which may be viewed as a \emph{generalized equilibrium set} evaluated at time $t\curr$. On the other hand, rather than finding closure to the discrete system~\eqref{eq:epsilon_disc}--\eqref{eq:eqmu_disc} by invoking constitutive relations (section~\ref{sec:mmorph_material}), we assume that the material response is encoded in a \emph{material data set} $\setD=\setD_1\times\dots\times\setD_M$, where $\setD_e\subset\mathrm{Z}_e$ is the local data set containing generalized stress-strain coordinates accessible to the material point $e$. We account for  history-dependent, inelastic responses by defining a suitable subset $\setD\curr\subset\setD$ that contains, at time $t\curr$, data points compatible with certain constraints imposed on the stress-strain history~\cite{eggersmann2019}. In this context, let us assume that a generic internal variable $\mathbf{q}$ is appended to the data set, such that a realization $\mathbf{q}_{\,e}$ is attached to every data point $\mathbf{z}_e\in\mathrm{D}_e$. Then, we consider a time-evolving data set evaluated at time $t\curr$:	
	\begin{equation}
			\setD_e\curr \coloneqq \Big\{ \mathbf{z}_e \in \setD_e \ \colon \  \Phi\big(\mathbf{q}_{\,e}\curr,\{(\eps_e,\gam_e,\zet_e,\sig_e,\btau_e,\bmu_e)^s\}_{s\leq k+1}\big)\leq \bm{0}  \Big\} \,,
	\label{eq:discdata_0}
	\end{equation}
where $\Phi$ is a generic function establishing a constraint to be fulfilled by the generalized stress-strain history and the internal variable state. Then, an evolving global data set can be defined at time $t\curr$ as $\setD\curr=\setD_1\curr\times\dots\times\setD_M\curr$. It bears emphasis that the availability of an internal state variable $\mathbf{q}_{\,e}$ does not entail the existence of a constitutive model; this variable must be measurable and included in the database alongside $\mathbf{z}_e\in\mathrm{D}_e$. \tR{For instance, in recent works on data-driven computing for granular materials, LS-DEM simulations~\cite{kawamoto2016} have been used by~\citet{karapiperis2021a} to obtain energy quantities and enforce thermodynamic constraints, and by~\citet{ulloa2023} to measure a breakage state variable and enforce irreversibility conditions in crushable media}.
								
		With the above definitions in hand, the time-discrete data-driven problem may be written as
		\begin{equation}
				\inf_{\displaystyle\mathbf{y}\in\setD\curr}\,\inf_{\displaystyle\mathbf{z}\in\setE\curr} \Vert \mathbf{y}-\mathbf{z}\Vert^2 = \inf_{\displaystyle\mathbf{z}\in\setE\curr}\,\inf_{\displaystyle\mathbf{y}\in\setD\curr} \Vert \mathbf{y}-\mathbf{z}\Vert^2 \,.
		\label{eq:DDproblem}
		\end{equation}
This problem is typically solved using a fixed-point iteration scheme. Fixing an estimate for the mechanically admissible state $\mathbf{z}^*\in\setE\curr$, the problem
		\begin{equation}
				\mathbf{y}^*=\arginf_{\displaystyle\mathbf{y}\in\setD\curr}\,\Vert \mathbf{y}-\mathbf{z}^*\Vert^2
		\label{eq:DDproblem_y}
		\end{equation}
is solved using a local nearest neighbor search for the data points $\mathbf{y}_e\in\setD_e\curr$. Then, fixing the estimate $\mathbf{y}^*\in\mathrm{D}\curr$, the problem
		\begin{equation}
				\mathbf{z}^*=\arginf_{\displaystyle\mathbf{z}\in\setE\curr}\,\Vert \mathbf{y}^*-\mathbf{z}\Vert^2
		\label{eq:DDproblem_z}
		\end{equation}
is solved enforcing the mechanical constrains embedded in~\eqref{eq:Eset}. The kinematics are enforced directly by replacing the generalized strains~\eqref{eq:epsilon_disc}--\eqref{eq:zeta_disc} in the local metric $\Vert\mathbf{z}_e-\mathbf{y}^*_e\Vert$, using equations \eqref{eq:locmet}--\eqref{eq:eqmu_disc}. On the other hand, the balance equations~\eqref{eq:eqsig_disc}--\eqref{eq:eqmu_disc} are enforced by means of Lagrange multipliers collected in vectors $\discrete{\bm{\lambda}}^u$ and $\discrete{\bm{\lambda}}^\rchi$. In light of equation~\eqref{eq:globmet}, we obtain the Euler-Lagrange equations
		\begin{align}
		&\sum_{e=1}^M w_e\Big[\Bet\,\Ctens_e\,\big(\Be\,\discrete{\mathbf{u}}-\eps_e^*\big)+\Bgt\,\Dtens_e\,\big(\Bg\,\discrete{\mathbf{u}}-\mathbf{B}^{\rchi}_e\,\discrete{\bm{\rchi}}-\gam^*_e\big)\Big]=\bm{0}\,,\\
		&\sum_{e=1}^M w_e\Big[\Bzt\,\Atens_e\,\big(\Bz\,\discrete{\bm{\rchi}}-\zet^*_e\big)-\Nct\,\Dtens_e\,\big(\Bg\,\discrete{\mathbf{u}}-\mathbf{N}^{\rchi}_e\,\discrete{\bm{\rchi}}-\gam_e^*\big)\Big]=\bm{0}\,,\\
		&\sum_{e=1}^M w_e\Big[\Bet\,\sig_e + \Bgt\,\btau_e\Big]-\discrete{\mathbf{f}}^\mathrm{ext}=\bm{0}\,,\\		
		&\sum_{e=1}^M w_e\Big[\Bzt\,\bmu_e - \Nct\,\btau_e\Big]-\discrete{\mathbf{m}}^\mathrm{ext}=\bm{0}\,,\\
		&w_e\Ctens_e^{-1}\big(\sig_e-\sig^*_e\big)-w_e\Be\,\discrete{\bm{\lambda}}^u=\bm{0}\,,\\
		&w_e\Dtens_e^{-1}\big(\btau_e-\btau^*_e\big)-w_e\Bg\,\discrete{\bm{\lambda}}^u+w_e\Nc\,\discrete{\bm{\lambda}}^\rchi=\bm{0}\,,\\
		&w_e\Atens_e^{-1}\big(\bmu_e-\bmu^*_e\big)-w_e\Bz\,\discrete{\bm{\lambda}}^\rchi=\bm{0}\,.					
		\end{align}						
		Upon substitution and rearrangement, we obtain the following linear system to be solved for $(\discrete{\mathbf{u}},\discrete{\bm{\rchi}})$:
		\begin{align}	
		&\sum_{e=1}^M w_e\Big[\Bet\Ctens_e\Be + \Bgt\Dtens_e\Bg\Big]\discrete{\mathbf{u}} - \sum_{e=1}^M w_e\Big[\Bgt\Dtens_e\Nc\Big]\discrete{\bm{\rchi}}=\sum_{e=1}^M w_e\Big[\Bet\Ctens_e\,\eps^*_e + \Bgt\Dtens_e\gam^*_e\Big],\label{eq:uprob1}\\
		&\sum_{e=1}^M w_e\Big[\Nct\Dtens_e\Bg\Big]\discrete{\mathbf{u}} - \sum_{e=1}^M w_e\Big[\Bzt\Atens_e\Bz + \Nct\Dtens_e\Nc\Big]\discrete{\bm{\rchi}}=\sum_{e=1}^M w_e\Big[\Nct\Dtens_e\gam^*_e-\Bzt\Atens_e\,\zet^*_e\Big],
		\label{eq:uprob2}
		\end{align}		
		and the following linear system to be solved for $(\discrete{\bm{\lambda}}^u,\discrete{\bm{\lambda}}^\rchi)$:
		\begin{align}	
		&\sum_{e=1}^M w_e\Big[\Bet\Ctens_e\Be + \Bgt\Dtens_e\Bg\Big]\discrete{\bm{\lambda}}^u - \sum_{e=1}^M w_e\Big[\Bgt\Dtens_e\Nc\Big]\discrete{\bm{\lambda}}^\rchi=\discrete{\mathbf{f}}^\mathrm{ext}-\sum_{e=1}^M w_e\Big[\Bet\sig^*_e + \Bgt\btau^*_e\Big],\label{eq:lprob1}\\
		& \sum_{e=1}^M w_e\Big[\Nct\Dtens_e\Bg\Big]\discrete{\bm{\lambda}}^u -\sum_{e=1}^M w_e\Big[\Bzt\Atens_e\Bz + \Nct\Dtens_e\Nc\Big]\discrete{\bm{\lambda}}^\rchi =\sum_{e=1}^M w_e\Big[\Bzt\bmu^*_e - \Nct\btau^*_e\Big] - \discrete{\mathbf{m}}^\mathrm{ext}.
		\label{eq:lprob2}		
		\end{align}
		The coupled systems~\eqref{eq:uprob1}--\eqref{eq:uprob2} and~\eqref{eq:lprob1}--\eqref{eq:lprob2} are standard forms for linear micromorphic elasticity and can be easily implemented in conventional finite element programs. The mechanical state is then updated~as
		\begin{align}
		&\begin{dcases}\eps_e=\Be\,\discrete{\mathbf{u}}\,,\\
		\sig_e=\Ctens_e\Be\,\discrete{\bm{\lambda}}^u + \sig_e^*\,;\end{dcases}\\
		&\begin{dcases}\gam_e=\Bg\,\discrete{\mathbf{u}}-\Nc\,\discrete{\bm{\rchi}}\,,\\
		\btau_e=\Dtens_e\Bg\,\discrete{\bm{\lambda}}^u -\Dtens_e\Nc\,\discrete{\bm{\lambda}}^\rchi + \btau_e^*\,;\end{dcases}\\
		&\begin{dcases}\zet_e=\Bz\,\discrete{\mathbf{\rchi}}\,,\\
		\bmu_e=\Atens_e\Bz\,\discrete{\bm{\lambda}}^\rchi + \bmu_e^*\,.\end{dcases}	
		\end{align}

This general formulation involves 12 (6) DOFs per node in 3D (2D), with $n=3$ ($n=2$) and $n_\chi=9$ ($n_\chi=4$), while the numerical operators $\Ctens_e$, $\Dtens_e$, and $\Atens_e$ involve 2, 3, and 11 parameters, respectively. However, as discussed in section~\ref{sec:mmorph_material}, simplifications of the general framework are possible. \tR{Hereafter, we consider, for the sake of simplicity, a microstrain continuum with $\bm\rchi$ purely symmetric (table~\ref{tab:mmorph}).} Then, the number of DOFs in 3D (2D) is reduced to 9 (5), with $n_\chi=6$ ($n_\chi=3$). Moreover, in agreement with remark~\ref{rem_1}, the total number of parameters in the operators $\Ctens_e$, $\Dtens_e$, and $\Atens_e$ is reduced to 4. \tR{Nevertheless, we emphasize that the full micromorphic version (or any subcase thereof) can also be implemented from the preceding formulation.}
    
\section{Numerical simulations}
\label{sec:num}

This section presents representative numerical simulations that showcase the ability of data-driven micromorphic computing to capture responses with strain localization. For the sake of conciseness, we limit the study to synthetic data sampled from the solution of reference BVPs. These data are then employed in data-driven computations that aim to replicate the reference model response. Hence, the main challenge for the data-driven algorithm will be to assign suitable data to each material~point.

In the examples that follow, a single data set is generated for all material points, such that $\mathrm{D}_e\curr\equiv\mathrm{D}\curr$, and the nearest neighbor search~\eqref{eq:DDproblem_y} is performed in a random order. Moreover, the history dependence of the material response is considered by imposing irreversibility conditions in the  data set~\eqref{eq:discdata_0}:
\begin{equation}
    \setD_e\curr \coloneqq \Big\{ \mathbf{z}_e \in \setD_e \ \colon \  0\leq q\curr-q\past \leq \mathtt{TOL}_q \Big\} \,.
    \label{eq:discdata_irr}
\end{equation}
Here, $q$ is an irreversible internal state variable, taken as the damage variable in sections~\ref{sec:bar} and~\ref{sec:V} and as the equivalent plastic strain in section~\ref{sec:biaxial}. Moreover, $\mathtt{TOL}_q$ is a tolerance parameter introduced to avoid large jumps in data sets that include loading-unloading and softening branches. 

\begin{remark}
    The data set parametrization~\eqref{eq:discdata_irr} assumes access to the internal state variable $q$, computed from the reference model in the current numerical study. However, for the application of this framework to real material data, this {irreversible} quantity must be accessible through empirical measurements. In this context, a measure of energy dissipation, computed from grain-scale simulations of granular materials, has been effectively employed~\cite{karapiperis2021a}. Additionally, a breakage state variable, measurable in experiments or grain-scale simulations, was used in a recent study~\cite{ulloa2023} focusing on crushable granular media.
\end{remark}

\subsection{One-dimensional softening bar with damage}
\label{sec:bar}

Let us first examine a simple 1D problem that allows us to identify the key issues related to softening responses and strain localization. In particular, we assess the standard data-driven framework for Cauchy continua in a softening bar, isolating the underlying challenges. These challenges are then addressed through the proposed extension to micromorphic continua.

\begin{figure}[t!]
    \centering 
    \includeinkscape[scale=0.95]{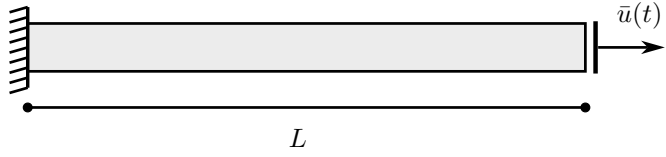}
    \caption{BVP for the 1D bar with constant cross-section.}
    \label{fig:bar1_bvp}
\end{figure}

Consider a 1D bar with constant unitary cross-section and length $L$ under tension, with applied displacements $u=\bar{u}(t)$ on one end and fixed displacements $u=0$ on the other (figure~\ref{fig:bar1_bvp}). The displacements are imposed in 153 increments from $\bar{u}=0$ to $\bar{u}=2.5$. The material is assumed to obey the gradient damage model with micromorphic elasticity presented in~\ref{sec:app1}, reduced to the 1D case ($n=1$ and $n_\chi=1$) through obvious model simplifications. The balance equations~\eqref{eq:dam1} and~\eqref{eq:dam2} are solved together with the damage evolution equations~\eqref{eq:dam3} using a standard staggered scheme. The constitutive parameters are chosen as $\mathsf{C}_{1111}/L=1$, $w_1L=1$, $\ell/L=0.8$, $c_1=1$, and $\ell_\rchi/L=0.1$. The bar is discretized into 100 finite elements using standard linear shape functions. Homogeneous Dirichlet boundary conditions are imposed on the damage field for the sake of obtaining a symmetric spatial response.
\begin{remark}
    In this normalized setting, arbitrary values of length $L$ will yield the same response of the model in terms of the nodal quantities $(u,L\rchi,\alpha)$ and state variables $(L\varepsilon, \sigma, L\gamma, \tau, L^2\zeta, \mu/L)$. For ease of notation, in this subsection, we let $\rchi\leftarrow L\rchi$, $\varepsilon\leftarrow L\varepsilon$, $\gamma\leftarrow L\gamma$, $\zeta\leftarrow L^2\zeta$, and $\mu\leftarrow \mu/L$.     
\end{remark}

\begin{figure}[b!]
    \centering 
    \includeinkscape[scale=0.85]{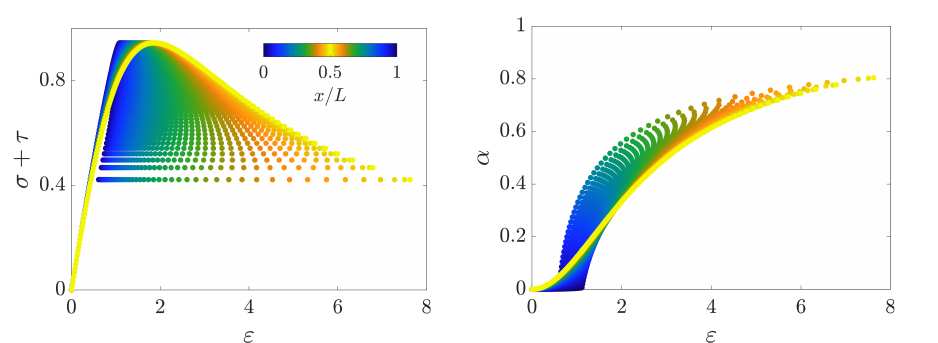}
    \caption{Data gathered from the solution of the 1D bar with constant cross-section, using the damage model in~\ref{sec:app1} simplified to the 1D case, to be used in standard data-driven computations. The data is taken from all the material points in the bar, yielding 15\,300 points (100 elements and 153 load steps). Note that $\sigma+\tau$ plays the role of Cauchy's~stress.}
    \label{fig:bar1_loc_data}    
\end{figure}

The standard data-driven framework for Cauchy continua is recovered from the general framework presented in section~\ref{sec:DD} by enforcing only the compatibility equation~\eqref{eq:epsilon_disc}
and the stress balance equation~\eqref{eq:eqsig}
in the equilibrium set~\eqref{eq:globmet}. The phase-space coordinates are then reduced to $\mathbf{z}_e=(\varepsilon_e,\sigma_e+\tau_e)$. As such, the stress quantity $\sigma_e+\tau_e$ is taken as an effective Cauchy stress, computed directly from the data to determine the material states and enforced to satisfy the standard equilibrium equations to determine the mechanical states. Accordingly, the distance metric~\eqref{eq:locmet}
takes the standard form 
\begin{equation}
 2\,\Vert \mathbf{z}_e \Vert^2 = E\varepsilon_e^2 + \frac{1}{E}(\sigma_e+\tau_e)^2\,,\
 \label{eq:met_1d_loc}
\end{equation}
where $E\coloneqq\mathsf{C}_{1111}$. Figure~\ref{fig:bar1_loc_data} presents the data points gathered from the solution of the BVP, showing softening in the central element with a transition to elastic unloading at $x=0$ and $x=L$. \tR{The data is taken from all the material points in the bar, yielding 15\,300 points (100 elements and 153 load steps).} The irreversible damage variable $\alpha$ is only employed to parametrize the history-dependent data set~\eqref{eq:discdata_irr}.

\begin{figure}[t!]
    \centering
    \footnotesize
    \includeinkscape[scale=0.85]{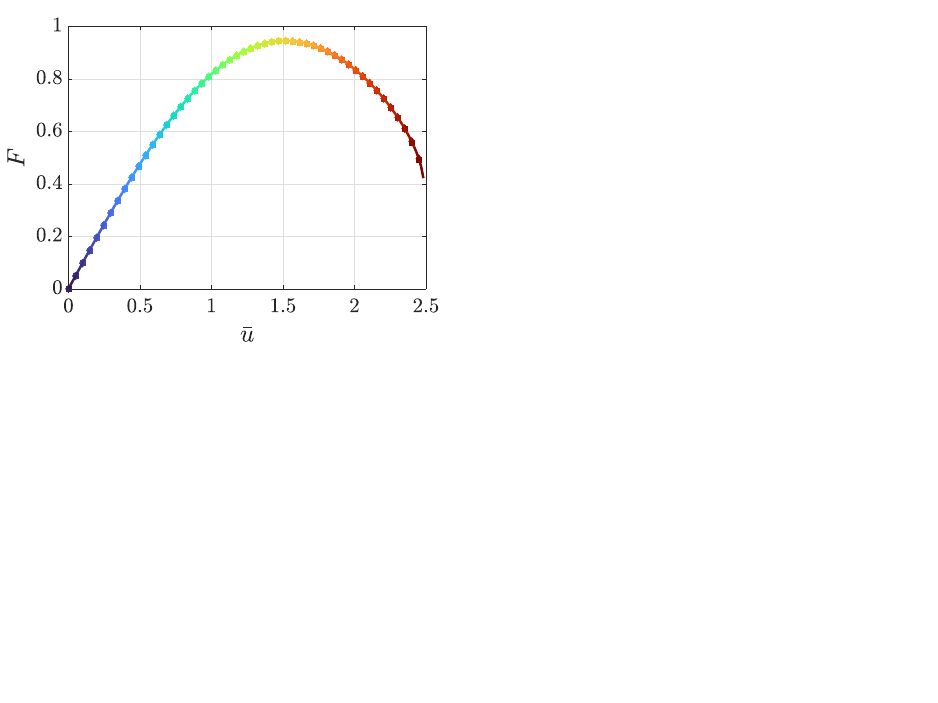}
    \caption{Data-driven results using the standard framework and comparison with the reference simulation (solid lines) for the 1D bar with constant cross-section, showing (a) the force-displacement response and the spatial distribution of (b) damage, (c) displacements, and (d) strains. The strain profile shows both the mechanical and material states, which visually overlap.}
    \label{fig:bar1_loc}
\end{figure}

\begin{figure}[t!]
    \centering 
    \footnotesize
    \includeinkscape[scale=0.85]{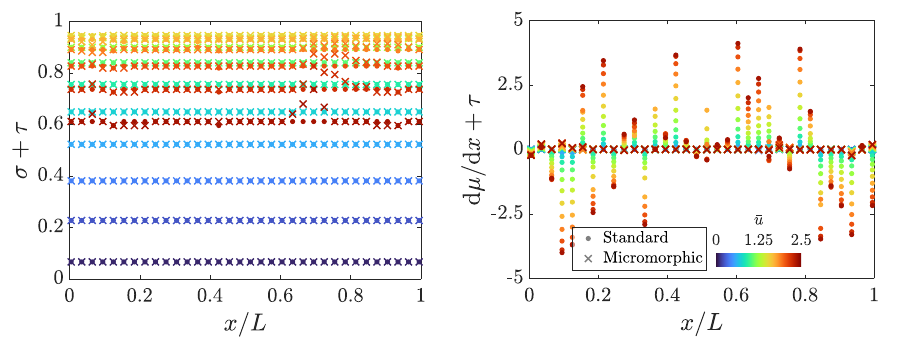}
    \caption{Verification of equilibrium conditions in the data-driven simulations of the 1D bar with constant cross-section, showing the results of the standard framework (dots) and the micromorphic framework (x marks). The color evolution corresponds to the load step $\bar{u}$. The spatial distribution of (a) $\sigma+\tau$ satisfies the equilibrium equation~\eqref{eq:eqsig_disc} with uniform profiles. The spatial distribution of (b) $\mathrm{d}\mu/\mathrm{d}\mathrm{x}+\tau$ satisfies the double-stress equilibrium equation~\eqref{eq:eqmu_disc} with vanishing profiles.}
        \label{fig:bar1_eq}
\end{figure}

Figure~\ref{fig:bar1_loc} shows the result of applying these data points in the standard data-driven framework for Cauchy continua. The initial undamaged parameter $E=C_{1111}=1$ is employed in the distance metric~\eqref{eq:met_1d_loc}. The global force-displacement curve shows a virtually perfect match between the reference response and the data-driven predictions, considering both the material state (i.e., the selected data points $\mathbf{y}^*\in\mathrm{D}\curr$) and the mechanical solution (i.e., $\mathbf{z}^*\in\mathrm{E}\curr$). On the other hand, the displacement, strain, and damage profiles (figures~\ref{fig:bar1_loc}b,c,d) are chaotic, in sharp contrast with the reference response. 

\begin{figure}[t!]
    \centering 
    \footnotesize
    \includeinkscape[scale=0.87]{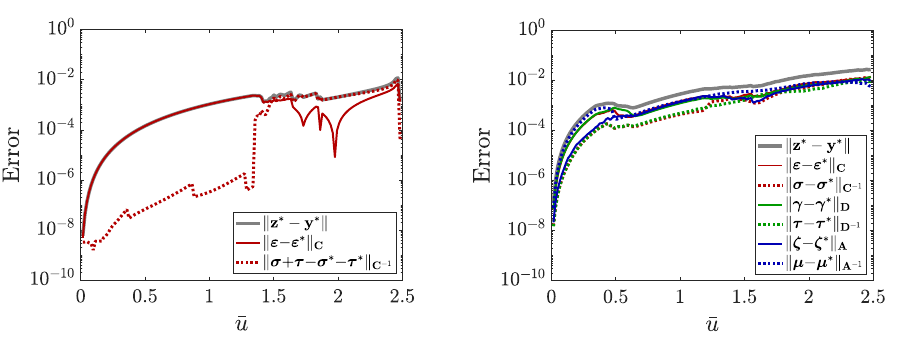}
    \caption{Evolution of the distance metric $\Vert\mathbf{z}^*-\mathbf{y}^*\Vert$ between material and mechanical states and individual contributions of the state variables for the data-driven simulations of the 1D bar with constant cross-section, showing the results of (a) the standard local framework and (b) the micromorphic framework.}
        \label{fig:bar1_dists}
\end{figure}

Note, nevertheless, that this solution is mechanically admissible, and corresponds to a local minimum of the strongly non-convex data-driven problem~\eqref{eq:DDproblem}. Notably, figure~\ref{fig:bar1_eq}a shows that the balance equation~\eqref{eq:eqsig_disc} is approximately satisfied at all time steps by the material state, up to the allowable error in figure~\ref{fig:bar1_dists}a. Indeed, any strain profile with constant stress values $\sigma_e+\tau_e$ along the length of the bar will comply with this condition. Therefore, different admissible profiles may emerge depending on the order in which the material points are read in the nearest neighbor search~\eqref{eq:DDproblem_y}. Then, unless we force the material points to only read the data corresponding to similar locations in the reference simulation, the standard algorithm has no reason to converge to a strain profile that captures the length scale of the reference response.

\begin{figure}[b!]
    \centering 
    \footnotesize
    \includeinkscape[scale=0.85]{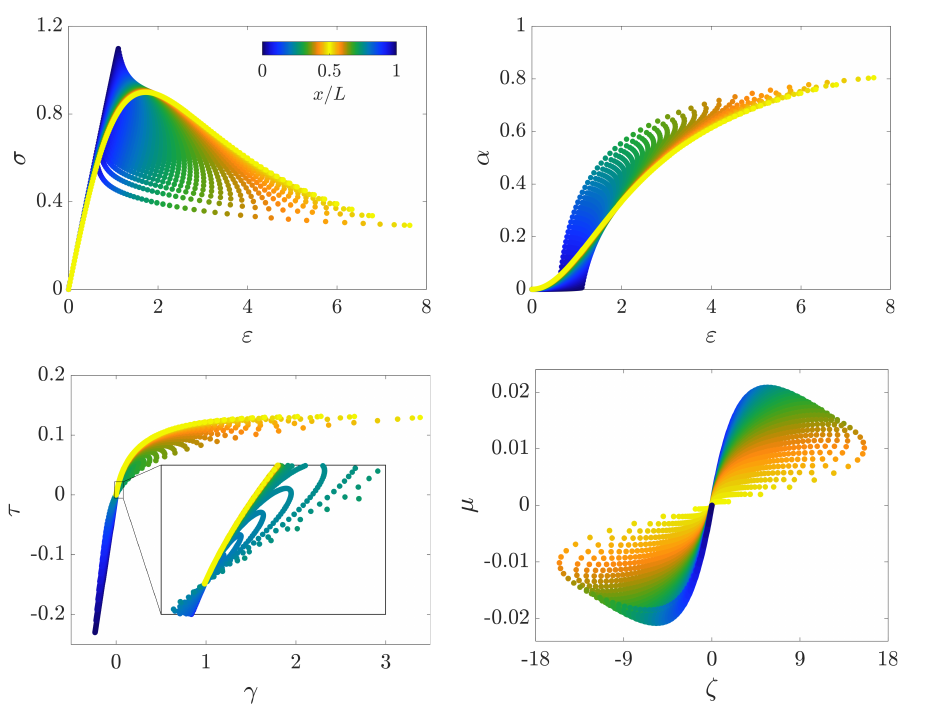}
    \caption{Data gathered from the solution of the 1D bar with constant cross-section, using the damage model in~\ref{sec:app1} simplified to the 1D case, to be used in data-driven micromorphic computations. The data is taken from all the material points in the bar, yielding a total of 15\,300 points (100 elements and 153 load steps).}
        \label{fig:bar1_mmorph_data}
\end{figure}

\begin{figure}[t!]
    \centering
    \footnotesize
    \includeinkscape[scale=0.85]{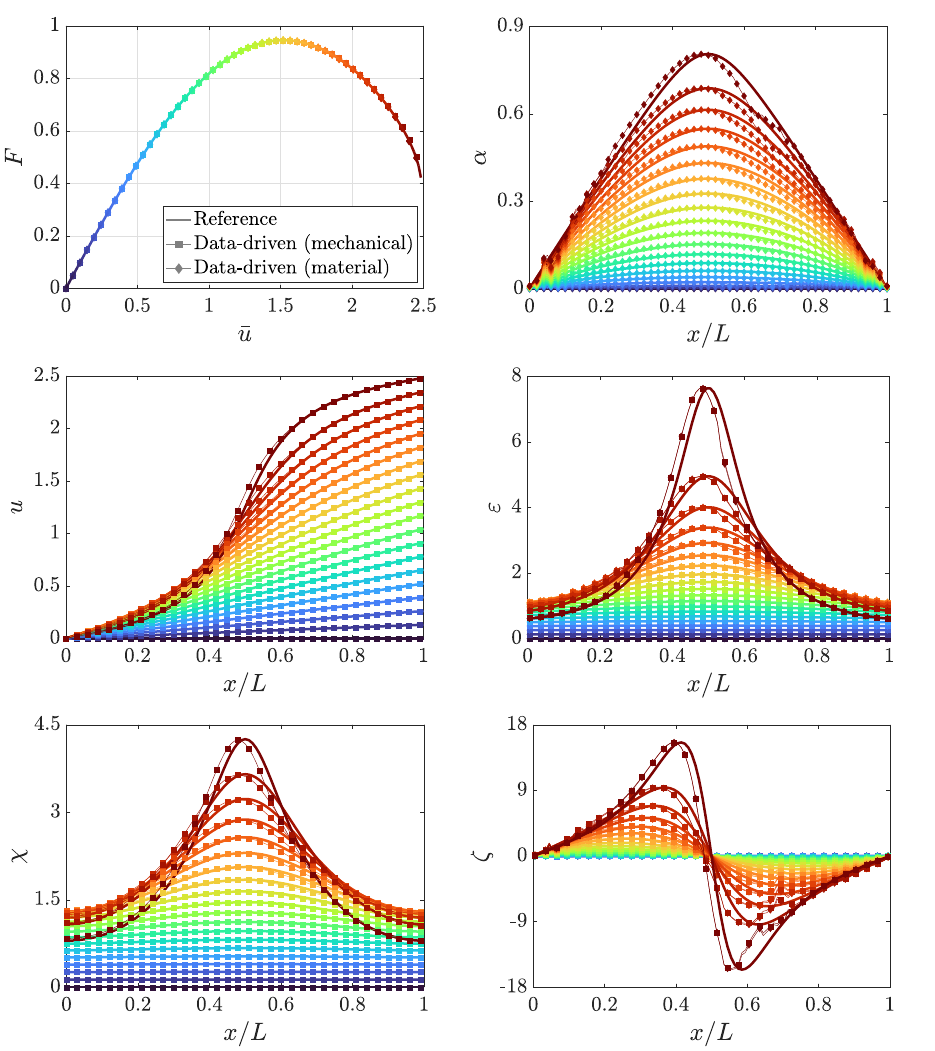}
    \caption{Data-driven results using the proposed micromorphic framework and comparison with the reference simulation (solid lines) for the 1D bar with constant cross-section, showing (a) the force-displacement response and the spatial distribution of (b) damage, (c) displacements, (d) strains, (e) micro-deformations, and (e) micro-deformation gradients. The strain and micro-deformation gradient profiles show both the mechanical and material states, which visually overlap.}
    \label{fig:bar1_mmorph}
\end{figure}

Let us now apply the proposed micromorphic formulation to this 1D problem. For this purpose, we proceed with the generalized phase-space coordinates $\mathbf{z}_e=(\varepsilon_e,\gamma_e,\zeta_e,\sigma_e,\tau_e,\mu_e)$. The full distance metric~\eqref{eq:locmet} then takes the form 
\begin{equation}
    2\,\Vert \mathbf{z}_e \Vert^2 = E\varepsilon_e^2 + \frac{1}{E}\sigma_e^2+c_1E\gamma_e^2 + \frac{1}{c_1E}\tau_e^2+c_1\ell_\rchi^2E\zeta_e^2 + \frac{1}{c_1\ell_\rchi^2E}\mu_e^2\,.
    \label{eq:met_1D_mmorph}
\end{equation}
Figure~\ref{fig:bar1_mmorph_data} presents the generalized data points collected from the reference solution. Again, the irreversible damage variable $\alpha$ is only employed to parametrize the history-dependent data set~\eqref{eq:discdata_irr}. As in figure~\ref{fig:bar1_loc_data}, the central element experiences softening, with a transition to elastic unloading at $x=0$ and $x=L$ in $(\varepsilon,\sigma)$ space. On the other hand, the central element exhibits a monotonically increasing behavior in $(\gamma,\tau)$ space and remains small in $(\zeta,\mu)$ space. For the latter, peak absolute values in $\zeta$ and $\mu$ occur at a distance from the center that evolves during the simulation, ending at $x/L\approx0.5\pm0.085$ and $x/L\approx0.5\pm0.25$, respectively. These profiles are better illustrated in the reference model response (figure~\ref{fig:bar1_mmorph}).

Figure~\ref{fig:bar1_mmorph} shows the result of applying these data points in the generalized data-driven framework for micromorphic continua. Again, the undamaged initial constitutive parameters are used in the distance metric~\eqref{eq:met_1D_mmorph}. As in figure~\ref{fig:bar1_loc}, the global force-displacement curve shows a very close agreement between the reference solution and both the material ($\mathbf{y}^*\in\mathrm{D}\curr$) and mechanical  ($\mathbf{z}^*\in\mathrm{E}\curr$) data-driven predictions. However, in this case, a realistic spatial distribution profile is predicted for all variables. In particular, figures~\ref{fig:bar1_mmorph}b--f show that the damage, displacement, strain, microstrain, and microstrain gradient profiles predicted by the data-driven algorithm are in close agreement with the reference response. 

Note, nevertheless, that a slight shift to the left is observed. As in the standard case, this result is a local minimum of the strongly non-convex data-driven problem~\eqref{eq:DDproblem}. However, the possible profiles are now restricted by the fulfillment of the double balance equation~\eqref{eq:eqmu_disc}. For illustration, Figure~\ref{fig:bar1_eq}b shows that the material state predicted by the micromorphic framework complies with equation~\eqref{eq:eqmu_disc}, up to the allowable error in figure~\ref{fig:bar1_dists}b, while the state predicted by the standard framework does not. Hence, while profiles with peaks at different locations are possible and observed in data-driven micromorphic simulations, the response must satisfy all the constraints embedded in the generalized mechanical manifold~\eqref{eq:Eset}, rendering profiles like those in figure~\ref{fig:bar1_loc} no longer admissible. We further note that the length scale of the reference response is now captured: while the stiffness encoded in the $(\varepsilon,\sigma)$ data is extracted by both frameworks, the length scale encoded in the $(\zeta,\mu)$ data is only extracted in the micromorphic simulations.

\begin{figure}[t!]
    \centering 
    \footnotesize
    \hspace{-2em}
    \includeinkscape[scale=0.87]{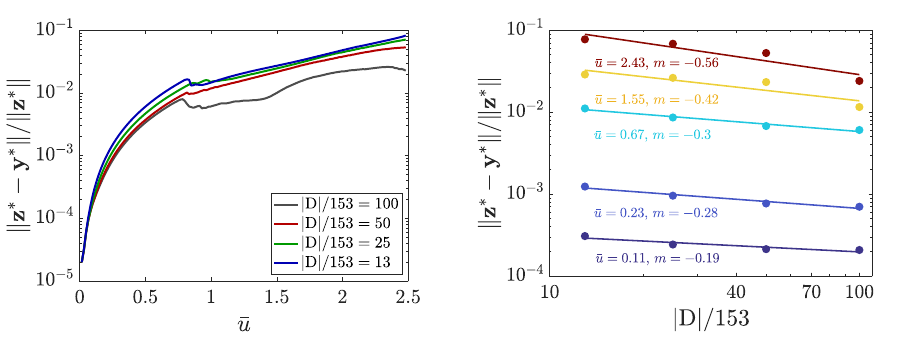}
    \caption{Relative distance between the material and mechanical states with increasing number of data points $|\mathrm{D}|$ in the data-driven micromorphic simulations of the 1D bar with constant cross-section, showing (a) the evolution of the relative distance metric and (b) the convergence rate in logarithmic space at different load steps $\bar{u}$. The cardinality $\vert\mathrm{D}\vert$ is divided by the number of load steps (153) to yield the number of sampling points.}
        \label{fig:bar1_conv}
\end{figure}

Figure~\ref{fig:bar1_conv} shows the convergence of the data-driven micromorphic simulations with respect to the number of data points. As expected, the relative error, measured as the relative distance between the material and mechanical states, i.e., $\Vert\mathbf{z}^*-\mathbf{y}^*\Vert/\Vert\mathbf{z}^*\Vert$, is smallest when the data is sampled from all the material points in the reference solution ($|\mathrm{D}|/153=100$). The data-driven predictions become less accurate as the applied displacement $\bar{u}$ increases due to the increasingly non-homogeneous response. Moreover, figure~\ref{fig:bar1_conv}b shows that the importance of the data size increases with $\bar{u}$. At the peak load $\bar{u}=1.5$, the relative error is approximately $2\%$ with only 13 sampled material points and approximately $0.9\%$ with all the data points. The corresponding error values increase to $8\%$ and $2\%$ at the final load~step.

\begin{figure}[b!]
    \centering 
    \includeinkscape[scale=0.95]{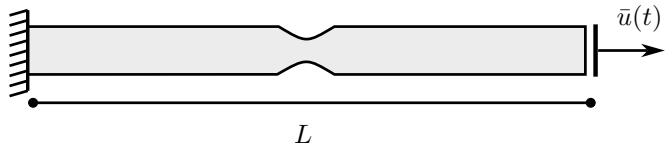}
    \caption{BVP for the 1D bar with a central notch. The cross-section is reduced to 0.5 in the middle along the coordinates $x/L\in[0.45,0.55]$.}
    \label{fig:bar2_bvp}
\end{figure}

\begin{figure}[t!]
    \centering
    \footnotesize
    \includeinkscape[scale=0.85]{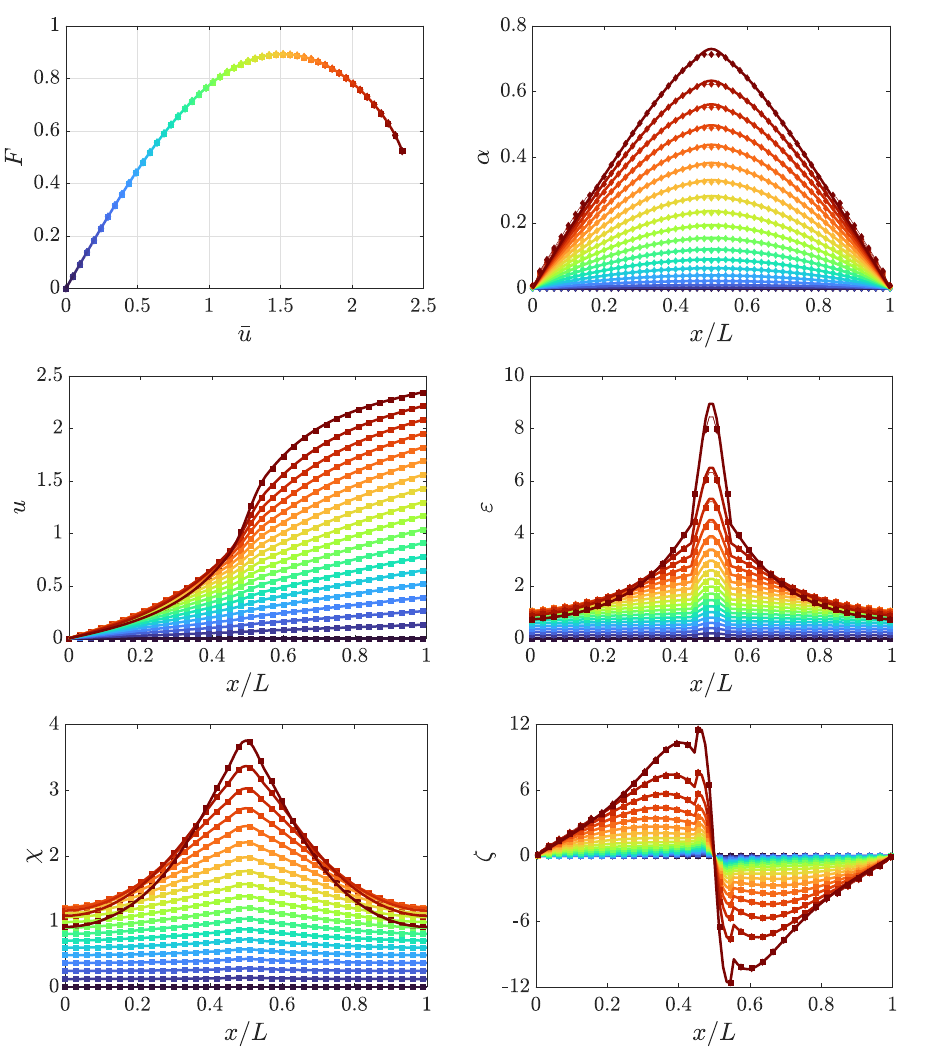}
    \caption{Data-driven results using the proposed micromorphic framework and comparison with the reference simulation (solid lines) for the 1D notched bar, showing (a) the force-displacement response and the spatial distribution of (b) damage, (c) displacements, (d) strains, (e) micro-deformations, and (f) micro-deformation gradients. The strain and micro-deformation gradient profiles show both the mechanical and material states, which visually overlap.}
    \label{fig:bar2_mmorph}
\end{figure}

Finally, we repeat the calculations above for a bar with a central geometrical notch (figure~\ref{fig:bar2_bvp}), imposing displacements in 145 increments from $\bar{u}=0$ to $\bar{u}=2.35$. Again, the standard framework fails to capture the spatial distribution of the reference response (these results are not shown for the sake of brevity). On the other hand, the data-driven prediction obtained with the micromorphic framework exhibits, in this case, peak strains at the center of the bar, coinciding with the reference solution (figure~\ref{fig:bar2_mmorph}). Note that the data is now taken from the reference solution of the notched bar (plots analogous to figures~\ref{fig:bar1_loc_data} and~\ref{fig:bar1_mmorph_data}, also not shown for the sake of brevity).

\begin{figure}[t!]
    \centering 
    \footnotesize
    \hspace{-2em}
    \includeinkscape[scale=0.87]{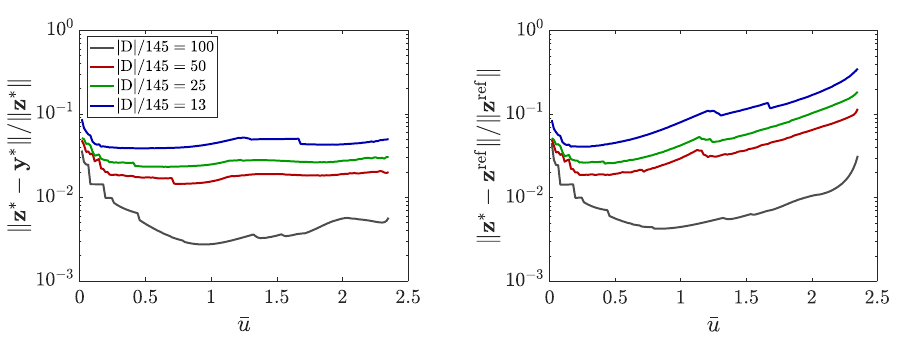}
    \caption{Relative distance with increasing number of data points $|\mathrm{D}|$ in the data-driven micromorphic simulations of the 1D notched bar, showing (a) the evolution of the relative distance metric between material and mechanical states, and (b) the evolution of the relative distance metric between the mechanical states and the reference solution $\mathbf{z}^\mathrm{ref}$.}
        \label{fig:bar2_conv}
\end{figure}

Figure~\ref{fig:bar2_conv}a shows the relative distance metric $\Vert\mathbf{z}^*-\mathbf{y}^*\Vert/\Vert\mathbf{z}^*\Vert$, quantifying the relative error between the material and mechanical states in the notched bar for different numbers of data points. At the end of the simulation, this value is approximately $5\%$ when only 13 material points are sampled and approximately $0.5\%$ when all 100 material points are sampled. Note that the error no longer increases monotonically as the load $\bar{u}$ increases. On the other hand, figure~\ref{fig:bar2_conv}b shows the relative errors with respect to the reference solution $\mathbf{z}^\mathrm{ref}$, i.e., the state extracted directly from the model response. At $\bar{u}=1.5$, the error ranges from $11.5\%$ with 13 sampled material points to $0.6\%$ with all the data points. At the final time step, these values increase, respectively, to approximately $35\%$ and $3\%$. Of course, because the relative distance between the mechanical and material states is much smaller, these responses are always mechanically admissible up to a certain tolerance.

At this point, we have shown that the proposed extension of data-driven mechanics to micromorphic continua can capture the non-homogeneous response of softening bars with no prior data assignment, whereas the (mechanically admissible) state predicted by the standard Cauchy continuum framework does not yield the expected profiles. The next sections evaluate the proposed framework in 2D scenarios.

\begin{remark}
 So far, the reference strain profiles have shown moderate localization, with none of the bars studied above reaching a damage level of 1 (i.e., the bars did not crack). While the proposed framework can handle scenarios with complete damage, as shown in the following example, this problem merits a word of caution. Specifically, the data set becomes sparse in these stages, with very low stresses and high strains. A dedicated approach may be considered in future works to properly navigate through such sparse data sets, e.g., performing adaptive sampling~\cite{gorgogianni2023}. 
 \label{rem:fulldam}
\end{remark}


\begin{figure}[b!]
    \centering 
    \vspace{1em}
    \small
    \includeinkscape[scale=0.74]{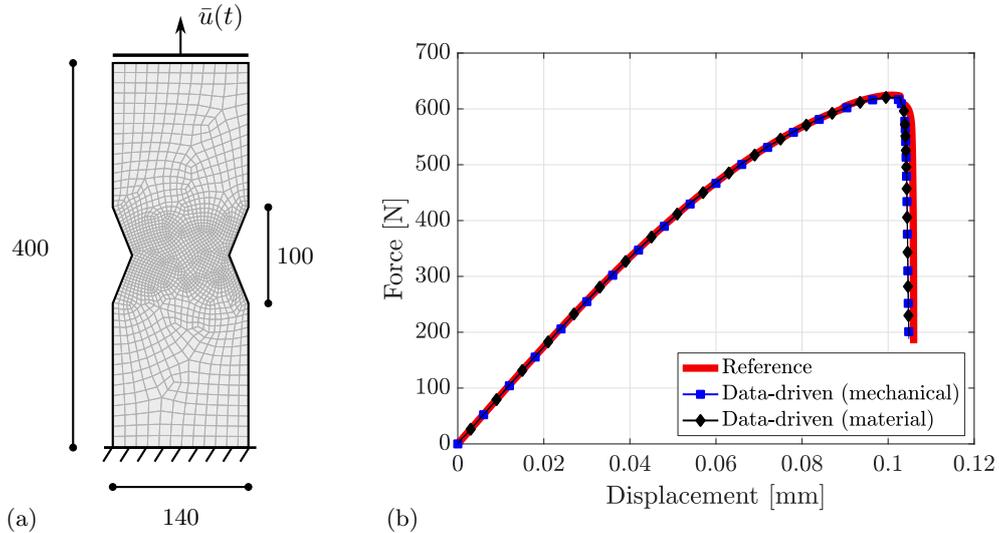}
	\put(-13cm,0.0cm){\tR{\small(a)}}
    \put(-8cm,0.0cm){\tR{\small(b)}}
    \vspace{1em}
    \caption{(a) BVP for the V-notched specimen with dimensions in mm and (b) force-displacement response comparing the data-driven solution to the reference simulation.}
    \label{fig:V_bvp_fd}
\end{figure}

\subsection{V-notched specimen with damage}
\label{sec:V}

We now apply the proposed data-driven micromorphic framework to a 2D BVP where the material experiences complete damage. In particular, consider a V-notched sample under plane-strain conditions and tensile loading (figure~\ref{fig:V_bvp_fd}a). Displacements are applied incrementally from $\bar{u}=0$ to $\bar{u}=0.106$~mm in 528 load steps. The material is assumed to obey the gradient damage model with micromorphic elasticity presented in~\ref{sec:app1}. We employ the following material parameters: Young's modulus $E=25.85$~GPa, Poisson's ratio $\nu=0.18$, damage threshold $w_1=4.8$~kPa, damage length scale $\ell=12.5$~mm, micromorphic constant $c_1=1$, and micromorphic length scale~$\ell_\rchi=12.5$~mm. The 2D specimen is discretized into 1664 quadrilateral elements, refining the region between the notches where fracture is expected, and using standard bilinear shape~functions. 

\begin{figure}[t!]
    \centering 
    \vspace{1em}
    \small
    \includeinkscape[width=0.97\linewidth]{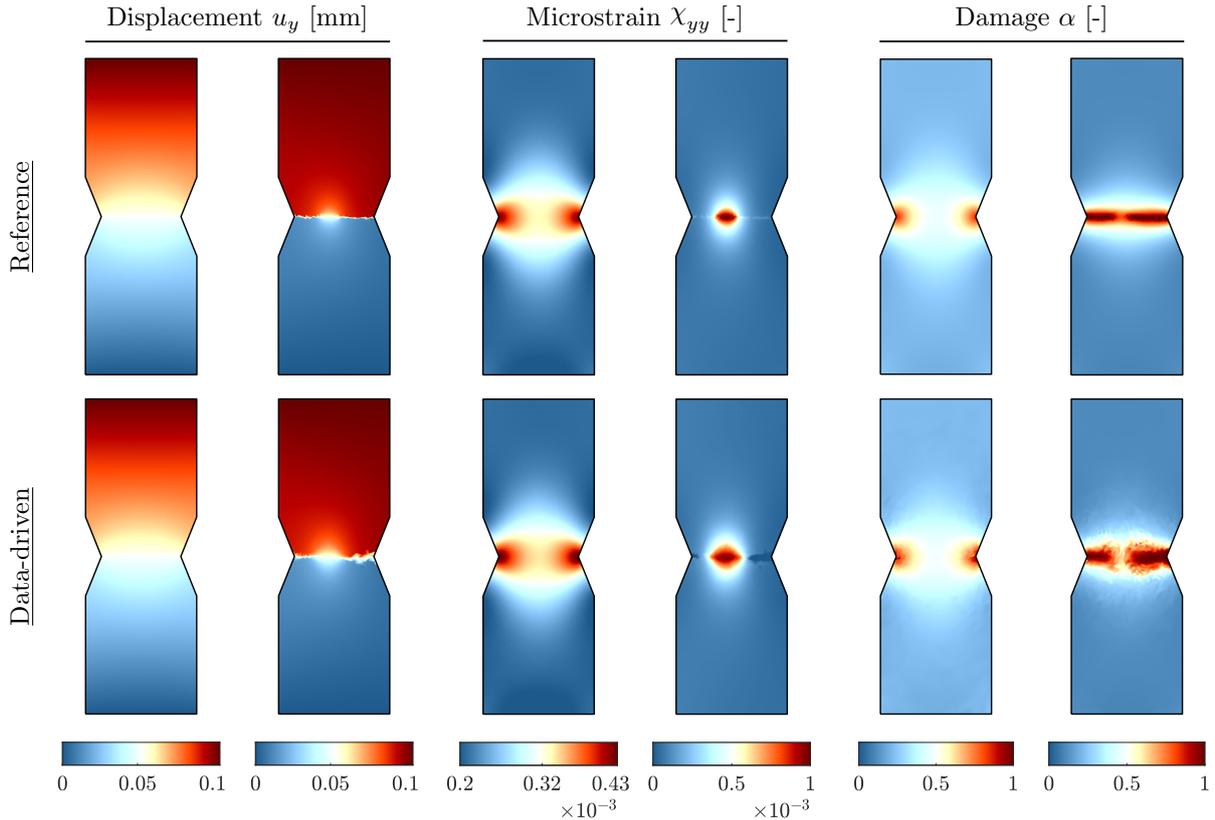}
    \caption{Reference simulation results for the V-notched specimen (top) and data-driven predictions (bottom), showing the vertical displacement component $u_y$ (left), the micro-deformation component  $\bm{\rchi}_{yy}$ (middle), and the damage field $\alpha$ (right). The snapshots correspond to initial crack propagation at the peak load and crack merging at failure.}
    \label{fig:V_kin}
\end{figure}

\begin{figure}[t!]
    \centering 
    \vspace{1em}
    \small
    \includeinkscape[width=0.97\linewidth]{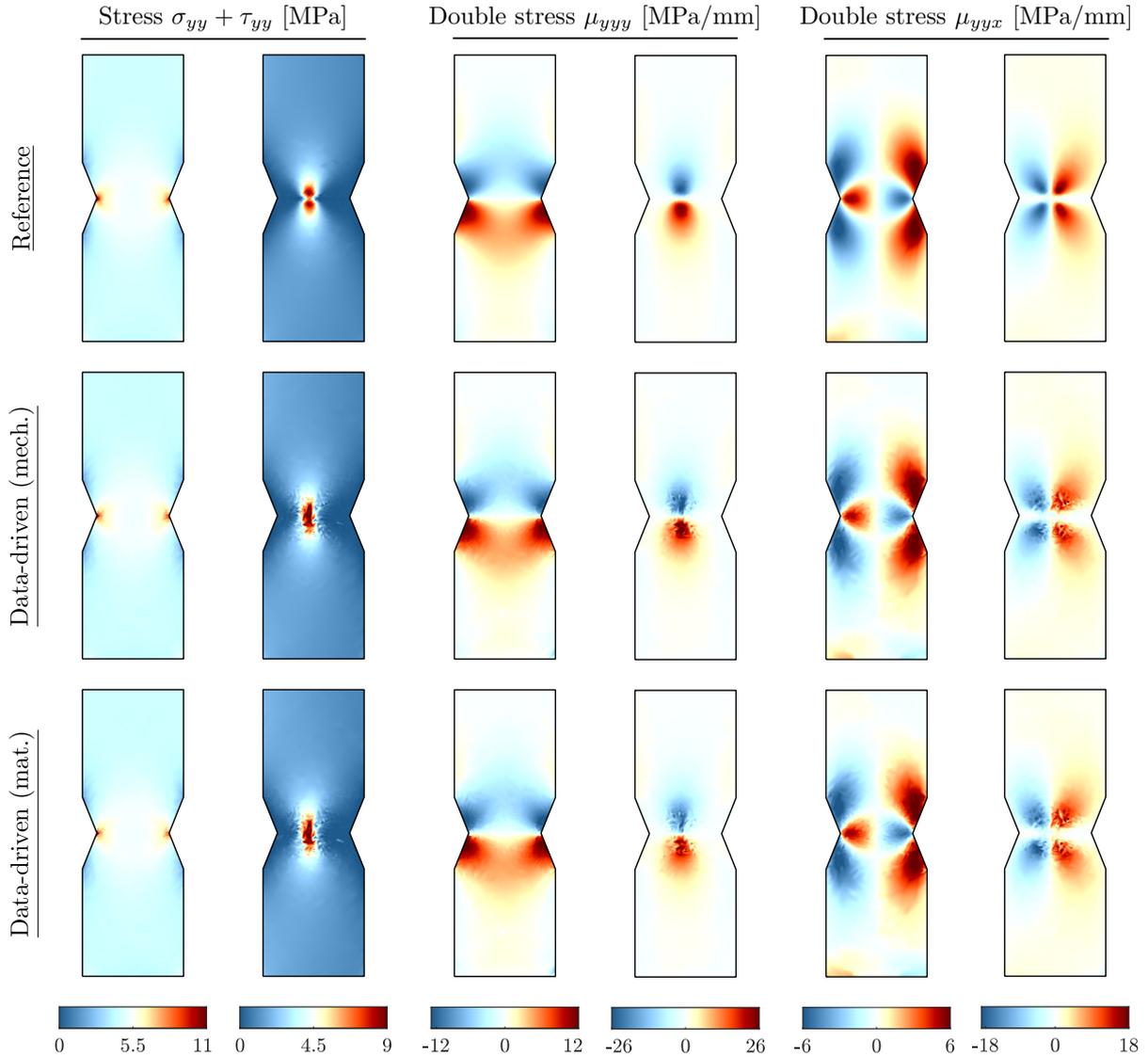}
    \vspace{1em}
    \caption{Reference simulation results for the V-notched specimen (top) and mechanical (middle) and material (bottom) data-driven predictions, showing the stress component $\sigma_{yy}+\tau_{yy}$ (left) and the double-stress components $\mu_{yyy}$ (middle) and $\mu_{yyx}$ (right). The snapshots correspond to initial crack propagation at the peak load and crack merging at failure.}
    \label{fig:V_str}
\end{figure}

The reference solution of the BVP is obtained by solving the governing equations~\eqref{eq:dam1}--\eqref{eq:dam3} of the micromorphic gradient damage model using a standard staggered scheme. Figure~\ref{fig:V_bvp_fd}b shows the force-displacement curve of the reference simulation (red line). The curve shows an initial, seemingly elastic stage. However, since the AT-2 formulation is used, with $w(\alpha)=w_1\alpha^2$, a slight level of damage evolves from the onset of loading and the response is never fully linear. The curve reaches a peak load at around $\bar{u}=0.1$~mm, followed by a softening response that leads to abrupt failure, characterizing brittle fracture. To reduce the brutal damage evolution at this stage, we include a slight numerical viscosity in the simulation, as commonly done in the phase-field literature~\cite{miehe2017}. Alternatively, dissipation-based arc-length techniques can be employed~\cite{wambacq2021}. Figure~\ref{fig:V_kin} (top) shows snapshots of the vertical displacement component $u_y$, the microstrain component $\rchi_{yy}$, and the damage variable $\alpha$. Cracks nucleate at the notches and propagate towards the middle of the sample. During this stage, the smoothness of the displacement field progressively decreases while the micro-deformation field concentrates at the crack tips. Figure~\ref{fig:V_str} (top) presents the corresponding stress fields, showing the $\sigma_{yy}+\tau_{yy}$, $\mu_{yyy}$, and $\mu_{yyx}$ components.

We proceed to reproduce the reference simulation using data-driven micromorphic computing. \tR{The generalized stress-strain states are extracted from the reference simulation at intervals of every 4 material points, i.e., utilizing 25\% of the available information. In total, we collect time-history data (load steps) from $1\,664$ sampling points, each yielding 528 generalized stress-strain states}. It is worth noting that the data becomes considerably sparse during the abrupt crack propagation stage, particularly at the final fracture steps when the specimen undergoes total failure. Addressing data sparsity in such responses is a challenge that merits consideration in future works. In this proof-of-concept simulation, we predict the response up to just before complete fracture, with a small ligament not fully broken. \tR{Additionally, we find it beneficial to update the metric tensors $\Ctens_e$, $\Dtens_e$, and $\Atens_e$ during this stage. In this simulation, we employ the initial undamaged moduli until the final propagation stage, where the secant moduli, which may be computed directly from the data, seem to yield better predictions. Alternatively, the use of optimal moduli~\cite{karapiperis2021b} may be considered.}

Figure~\ref{fig:V_bvp_fd}b compares the force-displacement response with the data-driven predictions, considering both the material ($\mathbf{y}^*\in\mathrm{D}\curr$) and mechanical  ($\mathbf{z}^*\in\mathrm{E}\curr$) states. The response shows a visually close agreement, even during the final crack propagation stage. Moreover, the proximity of the material and mechanical curves signals a relatively small distance metric throughout the simulation. Figures~\ref{fig:V_kin} and~\ref{fig:V_str} (bottom) show the corresponding spatial distribution of nodal quantities and generalized stress states. We observe a reasonable agreement between the reference response and the data-driven predictions. 

\begin{figure}[h!]
    \centering 
    \vspace{1em}
    \small
    \includeinkscape[width=0.97\linewidth]{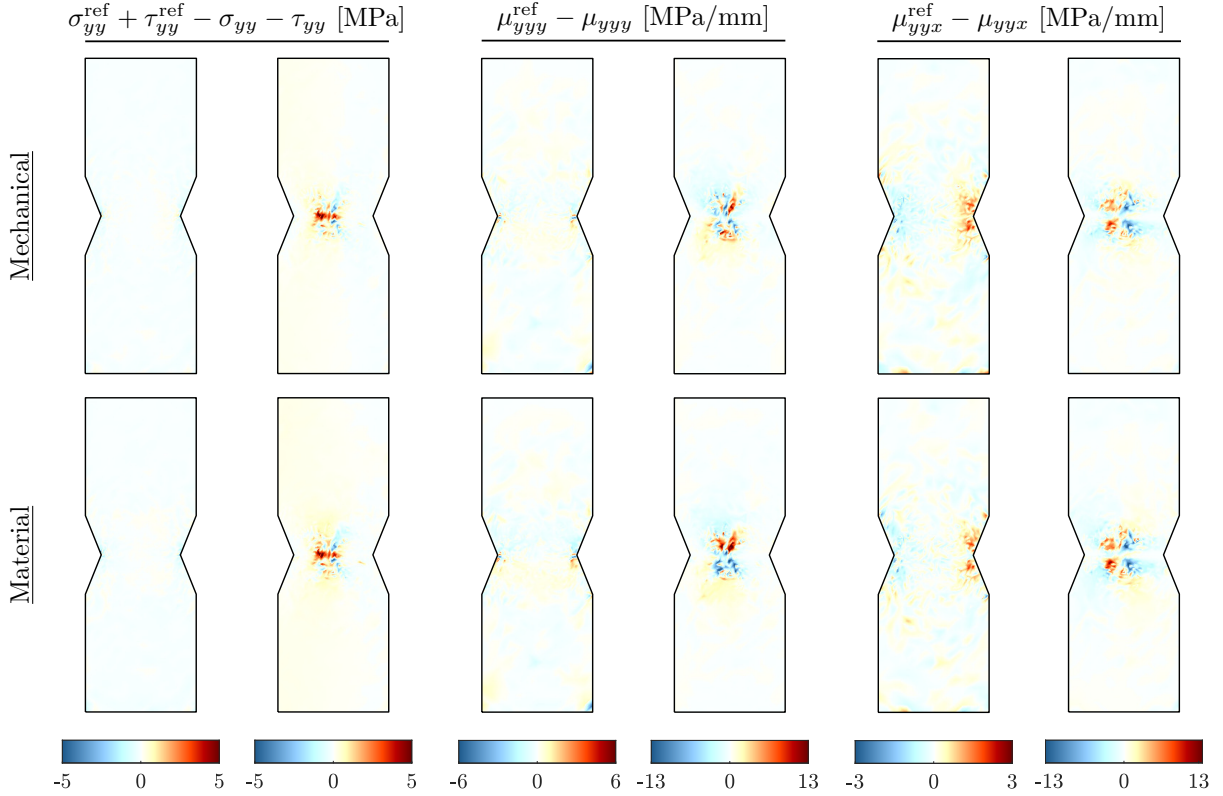}
    \caption{Deviation from the reference simulation for the V-notched specimen, showing the mechanical (top) and material (bottom) absolute error maps for the stress component $\sigma_{yy}+\tau_{yy}$ [MPa] (left) and the double-stress components $\mu_{yyy}$ [MPa/mm] (middle) and $\mu_{yyx}$ [MPa/mm] (right). The snapshots correspond to initial crack propagation at the peak load and crack merging at failure.}
    \label{fig:V_err}
\end{figure}

Note, however, that some noise is observed in the damage zone for the data-driven simulations, particularly at the final load step. Nevertheless, recall that the entire data set, consisting of over 1 million generalized stress-strain states, was read by all the material points with no particular assignment and in no particular order. Therefore, the deviations, mostly limited to specific zones near fractured elements, are deemed reasonable. Indeed, figure~\ref{fig:V_err} presents the absolute error in generalized stresses between the data-driven predictions and the reference simulation. The absolute error maps highlight regions of agreement and discrepancy in stress and double-stress components, with most errors limited to regions near the crack~tips. 

\subsection{Biaxial compression test with plasticity}
\label{sec:biaxial}

\begin{figure}[h!]
    \centering 
    \vspace{1em}
    \small
    \includeinkscape[scale=0.70]{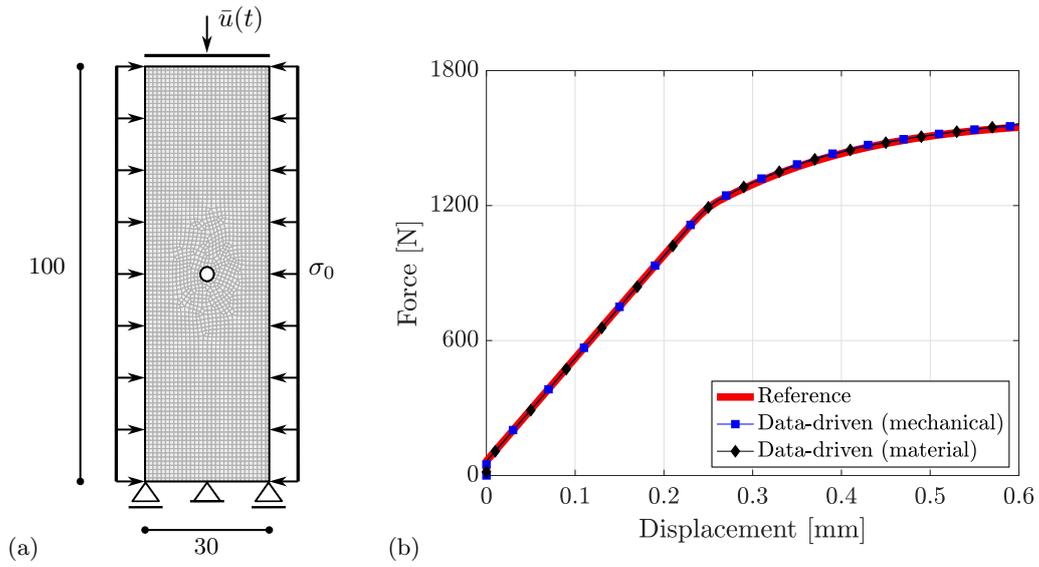}
	\put(-13.5cm,0.0cm){\tR{\small(a)}}
    \put(-8.5cm,0.0cm){\tR{\small(b)}}    
    \vspace{1em}
    \caption{(a) BVP for the biaxial compression test with dimensions in mm and (b) force-displacement response comparing the data-driven solution to the reference simulation.}
    \label{fig:biaxial_bvp_fd}
\end{figure}

The final example applies data-driven micromorphic computing to a  biaxial compression test under plane-strain conditions (figure~\ref{fig:biaxial_bvp_fd}), where the material is assumed to obey the micromorphic frictional-dilatant plasticity model presented in~\ref{sec:app2}. Vertical displacements are applied in 65 increments from $\bar{u}=0$ to $\bar{u}=0.6$ mm, keeping the lateral confining stress $\sigma_0=5$~MPa constant. We employ the following material parameters, taken from~\citet{borja2013} for the local version of the model: Young's modulus $E=14$~GPa, Poisson's ratio $\nu=0.3$, plastic yield strength $\sigma^\mathrm{p}=14.85$~MPa, friction coefficient $A_\varphi=0.29$, and dilation coefficient $A_\theta=0.25$. In addition, we consider a micromorphic constant $c_1=1$ and length scale~$\ell_\rchi=0.05$~mm. Finally, for the sake of simplicity, we take $a_1=a_2=a_3=a_4=0$.\footnote{Hence, only the classical plastic strain $\eps^\mathrm{p}$ is involved while ${\gam}^\mathrm{p}=\bm{0}$ and ${\zet}^\mathrm{p}=\bm{0}$. The same assumption is made in~\citet{dillard2006} in the context of metallic foams.} The 2D specimen is discretized into $3\,090$ quadrilateral elements using standard bilinear shape functions.

\begin{figure}[h!]
    \centering 
    \vspace{1em}
    \small
    \includeinkscape[width=0.95\linewidth]{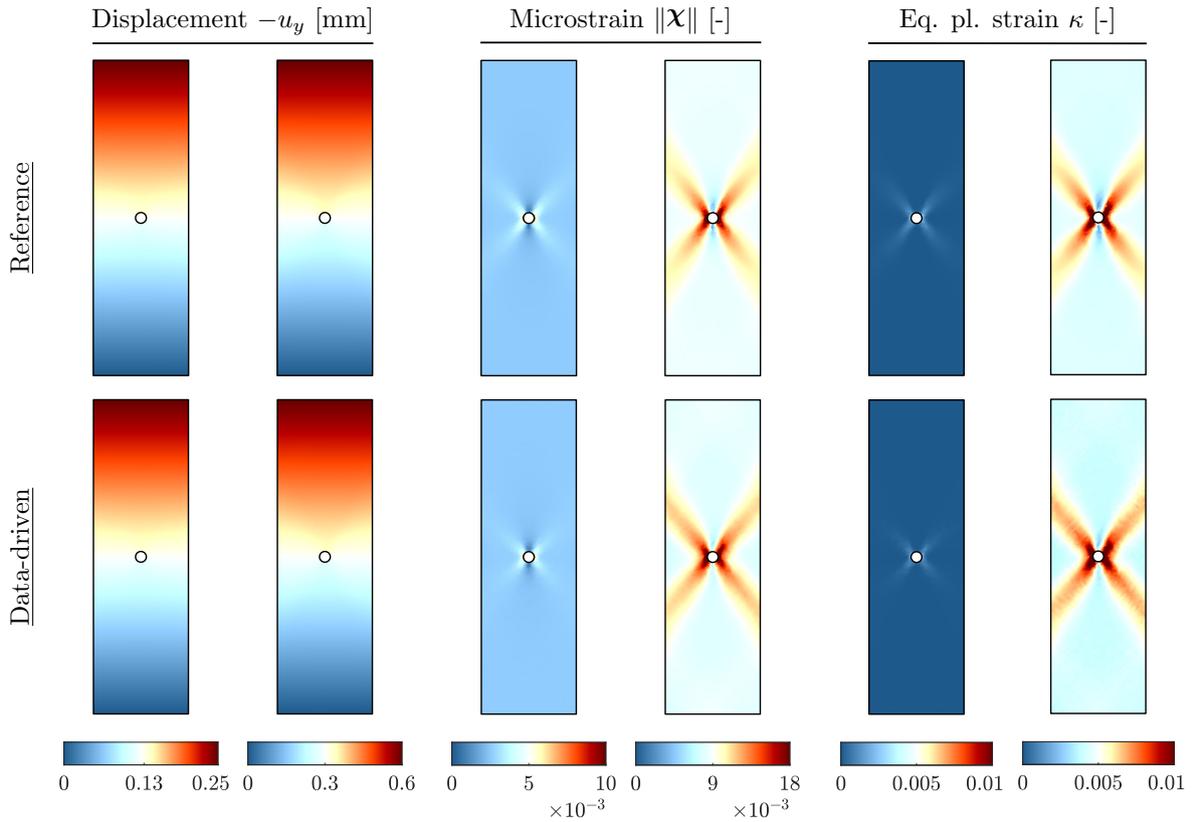}
    \vspace{1em}
    \caption{Reference simulation results for the biaxial compression test (top) and data-driven predictions (bottom), showing the vertical displacement component $-u_y$ [mm] (left), the micro-deformation $\Vert\bm{\rchi}\Vert$ (middle), and the equivalent plastic strain field $\kappa$ (right). The snapshots correspond to the initial yielding stage at $\bar{u}=0.25$ mm and the final load step at $\bar{u}=0.6$ mm.}    
    \label{fig:biaxial_kin}
\end{figure}

We proceed to reproduce the reference simulation using data-driven micromorphic computing, \tR{taking the time-history data from 25\% of all material points in a random order (65 time steps from 3\,090 sampling points)}. Figure~\ref{fig:biaxial_bvp_fd}b shows the global force-displacement curve measured on the top edge of the sample, comparing the reference solution to the material ($\mathbf{y}^*\in\mathrm{D}\curr$) and mechanical  ($\mathbf{z}^*\in\mathrm{E}\curr$) data-driven predictions. A reasonable agreement is observed throughout the simulation, with the material and mechanical states virtually overlapping. Accordingly, the displacements $u_y$, microstrains $\Vert\bm{\rchi}\Vert$, and equivalent plastic strains $\kappa$ in figure~\ref{fig:biaxial_kin} show that the data-driven predictions capture the reference shear-banding response. In particular,  at the initial yielding stage, where $\bar{u}\approx0.25$~mm, the micro-deformations show a slight concentration near the central hole in an inclined pattern, coinciding with the equivalent plastic strains. As the material points progressively yield, shear deformation bands develop symmetrically and propagate toward the right and left sides. The data-driven predictions capture this deformation profile accurately, with no prior data assignment and all material points reading a single data set.


\section{Conclusions and perspectives}
\label{sec:conclusion}

We have presented a data-driven formulation for generalized micromorphic continua. Its enhanced capabilities compared to the standard Cauchy-continuum framework have been illustrated in problems with strain localization without assigning data sets to specific material points a priori. The response of a 1D softening bar, assumed to be governed by a gradient damage model, has been studied systematically. It is shown numerically that local minima found by the standard data-driven framework may capture the global response but not the spatial distribution of the state variables expected from the reference data. In contrast, the proposed framework captures both the global response and the generalized stress-strain profiles, effectively extracting the material length scale through the solution of the data-driven BVP. The framework has been further validated in 2D for a V-notched sample undergoing fracture in tension and a biaxial sample with shear-banding, showing, in both cases, accurate predictions. \tR{These results suggest that data-driven micromorphic mechanics offers promising opportunities for capturing post-critical, inelastic behavior including both the global response and the failure mode. Moreover, the proposed framework is general, including any subcase of the full first-order micromorphic theory.}

\tR{The study paves the way for subsequent developments on various fronts. Firstly, practical applications of the present framework (beyond the proof-of-concept, synthetic studies considered here) require access to stress-strain data from microstructured material responses, without resorting to a micromorphic constitutive model. \ref{sec:app_data} provides several possibilities in this direction. Moreover, the high-dimensional space of generalized stress-strain states poses additional challenges stemming from the \emph{curse of dimensionality}. Hence, the present framework may benefit from numerical schemes with dimensionality reduction~\cite{he2021deep, ibanez2018manifold} adapted to the micromorphic case. Moreover, improved numerical strategies, e.g., using tensor voting~\cite{eggersmann2021b}, manifold learning~\cite{bahmani2022manifold}, or adaptive sampling~\cite{gorgogianni2023} may be instrumental for accuracy, particularly with limited material data. Finally, the present study does not delve into the data structure, using, at all times, a single data set for all material points. This approach serves as a proof of concept but is by no means efficient. More dedicated techniques that efficiently navigate through generalized stress-strain data may be envisaged, perhaps adopting a clustering approach~\cite{eggersmann2021} or quantum computing~\cite{xu2024quantum}.}

\section*{Acknowledgements}

JU and JA acknowledge support for this research provided by US ARO funding through the Multidisciplinary University Research Initiative (MURI) grant no. W911NF- 19- 1- 0245. LS and MO gratefully acknowledge the financial support of NExT ISite program of Nantes Universit\'e through International Research Project (IRP) iDDrEAM.

\appendix

\section{Reference models and synthetic data}
\label{sec:app} 

\subsection{Micromorphic gradient damage model}
\label{sec:app1} 

This appendix presents the reference continuum model used in the examples of sections~\ref{sec:bar}--\ref{sec:V}. The model is an extension of the conventional gradient-damage/phase-field approach to fracture~\cite{bourdin2000,marigo2016} to micromorphic elasticity. Thus, we introduce, in addition to the displacement field $\bm{u}$ and the micro-deformation field $\bm{\rchi}$, a damage variable $\alpha(\bm{x},t)\in[0,1]$, with 0 corresponding to an undamaged material point and 1 characterizing a completely damaged state. For the sake of brevity, we consider displacement boundary conditions on $\Gamma^u_{\mathrm{D}}$ and omit external double forces and applied micro-deformations.

We derive the model in a variational form in terms of an internal energy functional $\mathcal{E}$ and dissipation power functional $\mathcal{R}$:
\begin{equation}
\mathcal{E}(\bm{u},\bm{\rchi},\alpha)\coloneqq\int_\Omega\psi(\eps,\gam,\zet,\alpha,\nabla\alpha)\,\mathrm{d}\bm{x}\,,
\quad \mathcal{R}(\dot{\alpha};\alpha)\coloneqq\int_\Omega\phi(\dot{\alpha};\alpha)\,\mathrm{d}\bm{x}\,.
\end{equation}
The free energy density reads
\begin{equation}
\psi(\eps,\gam,\zet,\alpha,\nabla\alpha)\coloneqq\frac{1}{2}\Big[\hat{\Ctens}(\alpha)\eps:\eps+\hat{\Dtens}(\alpha)\big(\nabla\bm{u}-\bm{\rchi}\big):\big(\nabla\bm{u}-\bm{\rchi}\big)+\hat{\Atens}(\alpha)\nabla\bm{\rchi}\trip\nabla\bm{\rchi} + w_1\,\ell^2\nabla\alpha\cdot\nabla\alpha\Big]\,,
\end{equation}
yielding the generalized stresses
\begin{equation}
\sig=\frac{\partial\psi}{\partial\eps}=\hat{\Ctens}(\alpha)\eps\,, \qquad \btau=\frac{\partial\psi}{\partial\gam}=\hat{\Dtens}(\alpha)\gam\,, \qquad  \bmu=\frac{\partial\psi}{\partial\zet}=\hat{\Atens}(\alpha)\zet\,.
\label{eq:gstressdam}
\end{equation}
The dissipation potential is given by
\begin{equation}
\phi(\dot{\alpha};\alpha)\coloneqq\begin{dcases}w'(\alpha)\,\dot{\alpha} \quad &\text{if } \dot{\alpha}\geq 0 \,,\\
+\infty \quad &\text{otherwise}\,,\end{dcases} 
\label{eq:phi_d}
\end{equation}
with $w(\alpha)=w_1\alpha^2$ (AT-2 formulation). Here, $\ell$ is the damage length scale, while the parameter $w_1=w(1)$ represents the energy dissipated in a complete damage process for a homogeneous volume~element.

One can show that the energy balance equation
\begin{equation}
\frac{\mathrm{d}}{\mathrm{d}t}\mathcal{E}(\bm{u},\bm{\rchi},\alpha) + \mathcal{R}(\dot{\alpha};\alpha) - \int_{\Gamma^u_{\mathrm{N}}}\bar{\bm{t}}\cdot\dot{\bm{u}}\, \mathrm{d}S - \int_{\Gamma^u_{\mathrm{D}}}\bm{t}\cdot\dot{\bar{\bm{u}}}\, \mathrm{d}S =0
\label{eq:EB-1}
\end{equation}
and the first-order stability condition
\begin{equation}
\begin{aligned}
\delta\mathcal{E}(\bm{u},\bm{\rchi},\alpha)(\tilde{\bm{u}},\tilde{\bm{\rchi}},\tilde{\alpha}) +  \mathcal{R}(\tilde{\alpha};\alpha)- \int_{\Gamma^u_{\mathrm{N}}}\bar{\bm{t}}\cdot\tilde{\bm{u}}\, \mathrm{d}S \geq0 \quad \forall\,\{\tilde{\bm{u}},\tilde{\bm{\rchi}},\tilde{\alpha}\}\in\tilde{\EuScript{U}}\times\tilde{\EuScript{X}}\times\tilde{\EuScript{A}},
\label{eq:DS-1}
\end{aligned}
\end{equation}
where $\tilde{\EuScript{U}}\coloneqq\mathrm{H}^1_0(\Omega;\mathbb{R}^{n})$,
$\tilde{\EuScript{X}}\coloneqq\mathrm{H}^1_0(\Omega;\mathbb{R}^{n\times n})$, and 
$\tilde{\EuScript{A}}\coloneqq\mathrm{H}^1_0(\Omega;\mathbb{R}_+)$,
yield the following governing equations in strong form:
\begin{align}
&\text{Stress equilibrium} &&\begin{dcases}\dive\,(\sig+\btau)=\bm{0} \ \ \text{in} \ \ \Omega\,, \\
(\sig+\btau)\cdot\bm{n} = \bar{\bm{t}} \ \ \text{on} \ \ \Gamma^u_\mathrm{N}\,;
\end{dcases}\label{eq:dam1}\\
&\text{Double-stress equilibrium}
&&\begin{dcases}\dive\,\bmu+\btau=\bm{0} \ \ \text{in} \ \ \Omega\,,\\
\bm{\mu}\cdot\bm{n} = \bm{0}  \ \ \text{on} \ \ \Gamma\,;
\end{dcases}\label{eq:dam2}\\
&\text{Damage evolution}
&&\begin{dcases}-\frac{\partial\psi}{\partial\alpha}(\eps,\gam,\zet,\alpha,\nabla\alpha)-w'(\alpha)+w_1\ell^2\mathrm{div}[\nabla\alpha]\leq 0  \ \ \text{in} \ \ \Omega\,, \label{eq:dam3}\\ 
\bigg[-\frac{\partial\psi}{\partial\alpha}(\eps,\gam,\zet,\alpha,\nabla\alpha)-w'(\alpha)+w_1\ell^2\mathrm{div}[\nabla\alpha]\bigg]\dot{\alpha} = 0 \ \ \text{in} \ \ \Omega\,, \\ \dot{\alpha} \geq 0 \ \ \text{in} \ \ \Omega \,,\\ 
\nabla\alpha\cdot\bm{n} = {0}   \ \ \text{on} \ \ \Gamma\,.
\end{dcases}
\end{align}

For the sake of simplicity, following remark~\ref{rem_1}, we consider a microstrain formulation, with $\bm\rchi$ purely symmetric, and the following damaged elasticity tensors:
\begin{equation*}
    \hat{\mathsf{C}}_{ijkl}(\alpha)=(1-\alpha)^2\mathsf{C}_{ijkl}, \quad
\hat{\mathsf{D}}_{ijkl}(\alpha)=c_1\hat{\mathsf{C}}_{ijkl}(\alpha),  \quad \hat{\mathsf{A}}_{ijklmn}(\alpha)=c_1\ell_\rchi^2\hat{\mathsf{C}}_{ijlm}(\alpha)\delta_{kn}\,,
\end{equation*}
where $\mathsf{C}_{ijkl}$ is the standard isotopic elasticity tensor.

In sections~\ref{sec:bar}--\ref{sec:V}, we apply this model to solve different BVPs for $\bm{u}$ and $\bm{\rchi}$, from which the generalized stress/strain (synthetic) data is extracted as $$\begin{aligned}
\begin{dcases}
(\bm{\varepsilon},\bm{\sigma})=\big(\sym\nabla\bm{u},\,\hat{\Ctens}(\alpha)\sym\nabla\bm{u}\big)\,,\\
(\bm{\gamma},\bm{\tau})=\big(\nabla\bm{u}-\bm{\rchi},\,\hat{\Dtens}(\alpha)(\nabla\bm{u}-\bm{\rchi})\big)\,,\\
(\bm{\zeta},\bm{\mu})=\big(\nabla\bm{\rchi},\,\hat{\Atens}(\alpha)\nabla\bm{\rchi}\big)\,.
\end{dcases}
\end{aligned}$$

\subsection{Micromorphic plasticity model}
\label{sec:app2} 

This appendix presents the reference continuum model used in the examples of section~\ref{sec:biaxial}. The model is an extension of conventional non-associative  Drucker-Prager plasticity to the micromorphic setting. As in the previous section, we consider a microstrain continuum with $\bm\rchi$ purely symmetric. Moreover, similar to models presented in previous works~\cite{forest2004,forest2006,regueiro2009}, we may assume an additive decomposition of the micro-deformations into elastic and plastic components, such that $\bm{\rchi}=\bm{\rchi}^\mathrm{e}+\bm{\rchi}^\mathrm{p}$. The generalized strain measures are decomposed accordingly as
\begin{equation}
\eps=\eps^\mathrm{e}+\eps^\mathrm{p},\qquad\gam=\gam^\mathrm{e}+\gam^\mathrm{p},\qquad\zet=\zet^\mathrm{e}+\zet^\mathrm{p}\,.
\end{equation}
The free energy density reads
\begin{equation}
    \psi(\eps,\gam,\zet,\eps^\mathrm{p},\gam^\mathrm{p},\zet^\mathrm{p})\coloneqq\frac{1}{2}\Big[\Ctens\eps^\mathrm{e}:\eps^\mathrm{e}+\Dtens\gam^\mathrm{e}:\gam^\mathrm{e}+\Atens\zet^\mathrm{e}\trip\zet^\mathrm{e}\Big]\,.
\end{equation}
The conjugate stresses follow as
\begin{equation}
\sig^\mathrm{p}=-\frac{\partial\psi}{\partial\eps^\mathrm{p}}=\Ctens\eps^\mathrm{e}\equiv\sig\,, \qquad \btau^\mathrm{p}=-\frac{\partial\psi}{\partial\gam^\mathrm{p}}=\Dtens\gam^\mathrm{e}\equiv\btau\,, \qquad  \bmu^\mathrm{p}=-\frac{\partial\psi}{\partial\zet^\mathrm{p}}=\Atens\zet^\mathrm{e}\equiv\bmu\,.
\label{eq:gstressplas}
\end{equation}

At this point, different approaches may be taken to establish evolution laws for the three plastic internal variables $(\eps^\mathrm{p},\gam^\mathrm{p},\zet^\mathrm{p})$. In particular, as done in~\citet{forest2004} for compressible plasticity, we may consider either a combined yield function for all variables or an independent yield function for each variable. Focusing on non-associative Drucker-Prager plasticity, similar to~\citet{regueiro2009} and \citet{isbuga2017}, we may employ a combined yield function $f$ and plastic potential $g$:
\begin{equation}
    \begin{aligned}
    f(\sig,\btau,\zet)\coloneqq\sqrt{\frac{3}{2}\big(\Vert\mathrm{dev}\,\sig\Vert+a_1\Vert\mathrm{dev}\,\btau\Vert + a_3\Vert\mathrm{dev}\,\bmu\Vert\big)} + A_\varphi\big(\mathrm{tr}\,\sig+a_2\mathrm{tr}\,\btau+a_4\Vert\mathrm{tr}\,\bmu\Vert\big)-\sigma^\mathrm{p}\,,\\
    g(\sig,\btau,\zet)\coloneqq\sqrt{\frac{3}{2}\big(\Vert\mathrm{dev}\,\sig\Vert+a_1\Vert\mathrm{dev}\,\btau\Vert + a_3\Vert\mathrm{dev}\,\bmu\Vert\big)} + A_\theta\big(\mathrm{tr}\,\sig+a_2\mathrm{tr}\,\btau+a_4\Vert\mathrm{tr}\,\bmu\Vert\big)-\sigma^\mathrm{p}\,,
    \end{aligned}
    \label{eq:combyield}
\end{equation}
where $[\mathrm{dev}\,\bmu]_{ijk}=\mu_{ijk}-(1/3)\delta_{ij}\mu_{aak}$ and $[\mathrm{tr}\,\bmu]_{i}=\mu_{aai}$, or independent functions of the form
\begin{equation}
    \begin{aligned}
    &f^\varepsilon(\sig)\coloneqq\sqrt{\frac{3}{2}}\Vert\mathrm{dev}\,\sig\Vert + A_\varphi\mathrm{tr}\,\sig-\sigma^\mathrm{p}\,, \qquad &&g^\varepsilon(\sig)\coloneqq\sqrt{\frac{3}{2}}\Vert\mathrm{dev}\,\sig\Vert + A_\theta\mathrm{tr}\,\sig-\sigma^\mathrm{p}\,;\\
    &f^\gamma(\btau)\coloneqq\sqrt{\frac{3}{2}}a_1\Vert\mathrm{dev}\,\btau\Vert + a_2A_\varphi\mathrm{tr}\,\btau-a_5\sigma^\mathrm{p}\,, \qquad &&g^\gamma(\btau)\coloneqq\sqrt{\frac{3}{2}}a_1\Vert\mathrm{dev}\,\btau\Vert + a_2A_\theta\mathrm{tr}\,\btau-a_5\sigma^\mathrm{p}\,;\\
    &f^\zeta(\bmu)\coloneqq\sqrt{\frac{3}{2}}a_3\Vert\mathrm{dev}\,\bmu\Vert + a_4A_\varphi\Vert\mathrm{tr}\,\bmu\Vert-a_6\sigma^\mathrm{p}\,, \qquad &&g^\zeta(\bmu)\coloneqq\sqrt{\frac{3}{2}}a_3\Vert\mathrm{dev}\,\bmu\Vert + a_4A_\theta\Vert\mathrm{tr}\,\bmu\Vert-a_6\sigma^\mathrm{p}\,.
    \end{aligned}
    \label{eq:indepyield}
\end{equation}
Here, as in classical plasticity, $\sigma^\mathrm{p}\geq0$ is the yield strength, $A_\varphi\geq0$ is the friction coefficient, $A_\theta\in[0,A_\varphi]$ is the dilation coefficient reflecting the non-associativity of the model, and $a_1,\dots,a_6$ are additional material parameters for the micromorphic case.

A stress-dependent dissipation potential may now be evaluated following a generalization of the principle of maximum dissipation to non-associative models. For the case of independent functions~\eqref{eq:indepyield}: 
\begin{equation}
\phi(\dot{\eps}^\mathrm{p},\dot{\gam}^\mathrm{p},\dot{\zet}^\mathrm{p};\sig,\btau,\bmu)=\phi^\varepsilon(\dot{\eps}^\mathrm{p};\sig)+\phi^\gamma(\dot{\gam}^\mathrm{p};\btau)+\phi^\zeta(\dot{\zet}^\mathrm{p};\bmu),
\label{eq:phi_p}
\end{equation}
where (cf.~\cite{ulloa2021})
\begin{align}
\phi^\varepsilon(\dot{\eps}^\mathrm{p};\sig)&=\sup\big\{\tilde\sig:\dot\eps^\mathrm{p} \ \ \colon \ \ g^\varepsilon(\tilde\sig)\leq(A_\theta-A_\varphi)\tr\sig\big\}\\[0.5em]
&=\begin{dcases}\frac{\tr\dot{\eps}^\mathrm{p}}{3A_\theta}\big[(A_\theta-A_\varphi)\tr\sig+\sigma^\mathrm{p}\big] &\text{if} \quad \tr\dot{\eps}^\mathrm{p}\geq\sqrt{6}A_\theta\Vert\mathrm{dev}\,\dot{\eps}^\mathrm{p}\Vert\,, \\
+\infty &\text{otherwise};\end{dcases}\nonumber\\
\phi^\gamma(\dot{\gam}^\mathrm{p};\btau)&=\sup\big\{\tilde\btau:\dot\gam^\mathrm{p} \ \ \colon \ \ g^\gamma(\tilde\btau)\leq a_2(A_\theta-A_\varphi)\tr\btau\big\}\\[0.5em]
&=\begin{dcases}\frac{\tr\dot{\gam}^\mathrm{p}}{3a_2A_\theta}\big[a_2(A_\theta-A_\varphi)\tr\btau+a_5\sigma^\mathrm{p}\big] &\text{if} \quad \tr\dot{\gam}^\mathrm{p}\geq\frac{a_2}{a_1}\sqrt{6}A_\theta\Vert\mathrm{dev}\,\dot{\gam}^\mathrm{p}\Vert\,, \\
+\infty &\text{otherwise};\end{dcases}\nonumber\\
\phi^\zeta(\dot{\zet}^\mathrm{p};\bmu)&=\sup\big\{\tilde\bmu\trip\dot\zet^\mathrm{p} \ \ \colon \ \ g^\zeta(\tilde\bmu)\leq a_4(A_\theta-A_\varphi)\Vert\tr\bmu\Vert\big\}\\[0.5em]
&=\begin{dcases}\frac{\Vert\tr\dot{\zet}^\mathrm{p}\Vert}{3a_4A_\theta}\big[a_4(A_\theta-A_\varphi)\Vert\tr\bmu\Vert+a_6\sigma^\mathrm{p}\big] &\text{if} \quad \Vert\tr\dot{\zet}^\mathrm{p}\Vert\geq\frac{a_4}{a_3}\sqrt{6}A_\theta\Vert\mathrm{dev}\,\dot{\zet}^\mathrm{p}\Vert\,, \\
+\infty &\text{otherwise}.\end{dcases}\nonumber
\end{align}

With these functions in hand, one can show that the energy balance equation
\begin{equation}
\frac{\mathrm{d}}{\mathrm{d}t}\mathcal{E}(\bm{u},\bm{\rchi},\eps^\mathrm{p},\bm{\rchi}^\mathrm{p}) + \mathcal{R}(\dot{\eps}^\mathrm{p},\dot{\bm\rchi}^\mathrm{p};\sig) - \int_{\Gamma^u_{\mathrm{N}}}\bar{\bm{t}}\cdot\dot{\bm{u}}\, \mathrm{d}S - \int_{\Gamma^u_{\mathrm{D}}}\bm{t}\cdot\dot{\bar{\bm{u}}}\, \mathrm{d}S =0
\label{eq:EB-1_p}
\end{equation}
and the first-order stability condition
\begin{equation}
\begin{aligned}
\delta\mathcal{E}(\bm{u},\bm{\rchi},\eps^\mathrm{p},\bm{\rchi}^\mathrm{p})(\tilde{\bm{u}},\tilde{\bm{\rchi}},\tilde{\eps}^\mathrm{p},\tilde{\bm{\rchi}}^\mathrm{p}) +  \mathcal{R}(\tilde{\eps}^\mathrm{p},\tilde{\bm\rchi}^\mathrm{p};\sig,\btau,\bmu) - \int_{\Gamma^u_{\mathrm{N}}}\bar{\bm{t}}\cdot\tilde{\bm{u}}\, \mathrm{d}S \geq0 \quad \forall\,\{\tilde{\bm{u}},\tilde{\bm{\rchi}},\tilde{\eps}^\mathrm{p},\tilde{\bm{\rchi}}^\mathrm{p}\}\in\tilde{\EuScript{U}}\times\tilde{\EuScript{X}}\times\tilde{\EuScript{P}}\times\tilde{\EuScript{Y}},
\label{eq:DS-1_p}
\end{aligned}
\end{equation}
where
$\tilde{\EuScript{P}}\coloneqq\mathrm{L}^2(\Omega;\mathbb{R}^{n})$ and $\tilde{\EuScript{Y}}\coloneqq\mathrm{H}^1_0(\Omega;\mathbb{R}^{n\times n})$,
yield the following governing equations in strong form:
\begin{align}
&\text{Stress equilibrium} &&\begin{dcases}\dive\,(\sig+\btau)=\bm{0} \ \ \text{in} \ \ \Omega\,, \\
(\sig+\btau)\cdot\bm{n} = \bar{\bm{t}} \ \ \text{on} \ \ \Gamma^u_\mathrm{N}\,;
\end{dcases}\\
&\text{Double-stress equilibrium}
&&\begin{dcases}\dive\,\bmu+\btau=\bm{0} \ \ \text{in} \ \ \Omega\,,\\
\bm{\mu}\cdot\bm{n} = \bm{0}  \ \ \text{on} \ \ \Gamma\,;
\end{dcases}\\ 
&\text{Plasticity evolution}
&&\begin{dcases}f^\varepsilon(\sig)\leq 0, \quad f^\gamma(\btau)\leq 0, \quad f^\zeta(\bm{\mu})\leq 0  \ \ \text{in} \ \ \Omega\,, \\ 
\dot{\lambda}^\varepsilon\,f^\varepsilon(\sig)= 0, \quad \dot{\lambda}^\gamma\,f^\gamma(\btau)= 0, \quad \dot{\lambda}^\zeta\,f^\zeta(\bm{\mu})= 0  \ \ \text{in} \ \ \Omega\,, \\ 
\{\dot{\eps}^\mathrm{p},\dot{\gam}^\mathrm{p},\dot{\zet}^\mathrm{p}\}  \in\{\dot{\lambda}^\varepsilon \partial g^\varepsilon(\sig),\dot{\lambda}^\gamma \partial g^\gamma(\btau),\dot{\lambda}^\zeta \partial g^\zeta(\bm{\mu})\}  \ \ \text{in} \ \ \Omega\,.
\end{dcases}
\end{align}
The evolution equations for a model with the combined functions~\eqref{eq:combyield} can be derived similarly. Note further that a Cosserat-type plasticity model can be derived from the general formulation, taking $\bm\rchi$ purely skew-symmetric. Recent studies~\cite{de2022} have shown numerically that a model of this type resolves ill-posedness in BVPs caused by non-associativity.

{\color{black}

\section{\tR{On (generalized) stress-strain data acquisition}}\label{sec:app_data}

The apparent need for a generalized framework to capture responses with strain localization raises the natural question of whether the required generalized stress-strain data can be obtained from microstructured material responses in a practical scenario, i.e., without resorting to synthetic data from micromorphic constitutive models. We argue that homogenization towards generalized continua can be applied in concert with the present framework in applications to specific material systems. We discuss various possibilities below, compiling relevant references on the topic.
		 
		\subsection{From heterogeneous standard materials to generalized continua} 
		
		We first consider scenarios where the microstructure is described by standard, e.g., Cauchy-based, models, while the effective behavior is captured by generalized continua. Examples are heterogeneous materials with inclusions or voids and a weak separation of scales. For these cases, the link between the (classical and heterogeneous) microscopic behavior and the (generalized and homogeneous) macroscopic state has been established in several studies that can be taken as a point of departure.

		Depending on the application, we may pivot towards micropolar or second-gradient continua. In both cases, homogenization schemes are available to obtain generalized stress-strain states from the microscopic response. The homogenization of Cauchy-based heterogeneous materials towards micropolar continua has been addressed in several works~\cite{forest1998cosserat, feyel2003multilevel, de2011cosserat}, which provide effective micropolar stress-strain states from the response of representative volumes. The same applies to second-gradient continua, for which second-order homogenization~\cite{kouznetsova2002multi, kouznetsova2004multi, nguyen2013multiscale} schemes are available. In this case, we must deal with the requirement of $\mathrm{C}^1$ continuity in the numerical implementation.

		Of course, the present work addresses the more general case of micromorphic mechanics. In this broader setting, we may rely on several developments on the homogenization of heterogeneous materials towards micromorphic continua, some of which are listed below.
		
		\begin{itemize}
		\item An early contribution to micromorphic homogenization is found in the work of~\citet{forest2002homogenization}. Therein, kinematic averaging rules result from minimizing the difference between the actual, microscopic displacement $\bm{v}(\bm{x},\bm{\xi})$ and the approximating ansatz $\bm{u}(\bm{x})+\bm{\rchi}(\bm{x})\cdot\bm{\xi}$ (cf. figure~\ref{fig:micbvp}):
		\begin{equation}
		\min_{(\bm{u},\bm{\rchi})}\Big\langle \Vert \bm{v}(\bm{x},\bm{\xi}) - \bm{u}(\bm{x}) - \bm{\rchi}(\bm{x})\cdot\bm{\xi} \Vert^2 \Big\rangle_{\Omega^\mathrm{m}}\,,
		\label{eq:min0}
		\end{equation}
		leading to
		    \begin{equation}\label{eq:kin_hom0}
		       \bm{u}(\bm{x})=\big\langle\bm{v}(\bm{x},\bm{\xi})\big\rangle_{\Omega^\mathrm{m}}, \qquad         \bm{\rchi}(\bm{x}) = \big\langle\bm{v}(\bm{x},\bm{\xi})\otimes\bm{\xi}\big\rangle_{\Omega^\mathrm{m}}\cdot\big\langle\bm{\xi}\otimes\bm{\xi}\big\rangle^{-1}_{\Omega^\mathrm{m}}\,,
		    \end{equation}
		from which the gradients $\nabla_{\bm{x}}\bm{u}$ and $\nabla_{\bm{x}}\bm{\rchi}$ readily follow. The microdeformation is identified as the first moment of the microscopic displacement. An immediate issue is that these kinematic volume averages cannot be written as surface integrals, hindering the macro-micro transition, i.e., the possibility of applying the kinematic quantities as boundary conditions in representative volumes. This problem can be addressed by adopting polynomial ansatzes~\cite{forest2011, janicke2009two} based on certain ad hoc assumptions. A more consistent approach can be found in~\citet{janicke2012minimal}, employing the concept of \emph{minimal loading conditions}, where the averaging relations descending from~\eqref{eq:min0} (and the corresponding gradients) are imposed as global constraints in the microscale BVP. The corresponding generalized stresses follow as the Lagrange multipliers of the constrained BVP~\cite{hutter2019theory}. 
		
		\item While the frameworks above follow the micromorphic theory based on the principle of virtual power~\cite{mindlin1964,germain1973}, an alternative framework is found in~\citet{hutter2017homogenization,hutter2019theory}, which takes, instead, the averaging-based micromorphic formulation of~\citet{eringen1968mechanics} as the point of departure. In particular, noting that \citet{eringen1968mechanics} did not delve into the kinematics and did not provide explicit expressions for the average kinetics, \citet{hutter2017homogenization} derived consistently the following averaging rules:
		\begin{equation}
		\begin{aligned}
		&\bm{\sigma}(\bm{x}) = \big\langle\bm{\sigma}^\mathrm{m}(\bm{x},\bm{\xi})\big\rangle_{\Omega^\mathrm{m}}\,, \\
		&\bm{\tau}(\bm{x}) = \frac{1}{|\Omega^\mathrm{m}|}\int_{\partial\Omega^\mathrm{m}}\bm{\sigma}^\mathrm{m}(\bm{x},\bm{\xi})\cdot\bm{n}^\mathrm{m}(\bm{\xi})\otimes\bm{\xi}\,\mathrm{d}S - \big\langle\bm{\sigma}^\mathrm{m}(\bm{x},\bm{\xi})\big\rangle_{\Omega^\mathrm{m}}\,, \\
		&\bm{\mu}(\bm{x}) = \frac{1}{|\Omega^\mathrm{m}|}\int_{\partial\Omega^\mathrm{m}}\bm{\sigma}^\mathrm{m}(\bm{x},\bm{\xi})\cdot\bm{n}^\mathrm{m}(\bm{\xi})\otimes\bm{\xi}\otimes\bm{\xi}\,\mathrm{d}S\,,
		\end{aligned}
		\end{equation}
		with $\bm{\sigma}^\mathrm{m}$ the standard Cauchy stress tensor at the microscale and $\bm{n}^\mathrm{m}$ the unit normal at $\partial\Omega^\mathrm{m}$. The kinematic averaging rules are identical to~equation~\eqref{eq:kin_hom0} and the gradients thereof, to be imposed as \emph{global constraints}, adopting the concept of minimal loading conditions~\cite{janicke2012minimal}. 

		\item The previous approach is not suitable for materials with voids because it relies on kinematic volume constraints, requiring the microscopic displacement field to be defined everywhere inside the RVE. Indeed, in the microscale BVP, the generalized stresses appear as volume forces, which cannot be carried by voids. For such cases, \citet{hutter2019theory} proposed to replace the standard volume averaging $\langle\Box\rangle_{\Omega^\mathrm{m}}$ in equation~\eqref{eq:kin_hom0} with a modified kinematic averaging operator
		    \begin{equation}\label{eq:av_op}
		\big\langle\Box\big\rangle_\mathrm{M}\coloneq\frac{1}{|\Omega^\mathrm{m}|\big\langle H_\mathrm{M}(\bm{\xi})\big\rangle_{\Omega^\mathrm{m}}}\int_\mathrm{\Omega^\mathrm{m}} H_\mathrm{M}(\bm{\xi})\,\Box\,\mathrm{d}\bm{\xi}\,,
		    \end{equation}
		where $H_\mathrm{M}$ is a weighting function that vanishes inside the voids. In this way, a scheme fully compatible with FE$^2$ is established for various applications. 
		
		\item Another approach presented in~\citet{alavi2021construction} employs a quartic expansion for the microscopic displacement field. This method is able to establish surface integrals for all kinematic and kinetic quantities from variational principles.
		
		\item A different framework, suitable for composites with soft inclusions or voids, is found in~\citet{biswas2017micromorphic}. Therein, the microdeformation tensor is defined as the average distortion within the microheterogeneities, taking the form of a surface integral over internal interfaces:
		    \begin{equation}\label{eq:kin_hom1}
		      \bm{\rchi}(\bm{x}) = \big\langle\nabla_{\bm{\xi}}\bm{v}(\bm{x},\bm{\xi})\big\rangle_{\Omega^\mathrm{i}}=\frac{1}{|\Omega^\mathrm{i}|}\int_{\partial\Omega^\mathrm{i}}\bm{v}(\bm{x},\bm{\xi})\otimes\bm{n}^\mathrm{i}(\bm{\xi})\,\mathrm{d}S\,,
		    \end{equation} 
		 where $\Omega^\mathrm{i}\subset\Omega^\mathrm{m}$ is the domain covered by inclusions or voids. If this condition is imposed in the microscale BVP as
		 \begin{equation}
		\min_{(\bm{u},\bm{\rchi})}\bigg\{\int_{\Omega^\mathrm{m}} \sigma^\mathrm{m}(\bm{x},\bm{\xi}):\nabla_{\bm{\xi}}\bm{v}(\bm{x},\bm{\xi})\,\mathrm{d}\bm{\xi} - \bm{\lambda}(\bm{x}):\bigg(\bm{\rchi}(\bm{x})-\frac{1}{|\Omega^\mathrm{i}|}\int_{\partial\Omega^\mathrm{i}}\bm{v}(\bm{x},\bm{\xi})\otimes\bm{n}^\mathrm{i}(\bm{\xi})\,\mathrm{d}S\bigg)\bigg\}\,,
		\label{eq:min1}
		\end{equation}
		 together with suitable boundary conditions, the generalized stresses follow as 
		\begin{equation}
		\begin{aligned}
		&\bm{\sigma}(\bm{x}) = \frac{1}{|\Omega^\mathrm{m}|}\int_{\partial\Omega^\mathrm{m}}\bm{\sigma}^\mathrm{m}(\bm{x},\bm{\xi})\cdot\bm{n}^\mathrm{m}(\bm{\xi})\otimes\bm{\xi}\,\mathrm{d}S \, + \, \frac{|\Omega^\mathrm{i}|}{|\Omega^\mathrm{m}|}\bm{\lambda}(\bm{x})\,, \\
		&\bm{\tau}(\bm{x}) = -\frac{|\Omega^\mathrm{i}|}{|\Omega^\mathrm{m}|}\bm{\lambda}(\bm{x})\,, \\
		&\bm{\mu}(\bm{x}) = \frac{1}{2|\Omega^\mathrm{m}|}\int_{\partial\Omega^\mathrm{m}}\bm{\sigma}^\mathrm{m}(\bm{x},\bm{\xi})\cdot\bm{n}^\mathrm{m}(\bm{\xi})\otimes\bm{\xi}\otimes\bm{\xi}\,\mathrm{d}S\,,
		\end{aligned}
		\end{equation} 
		encompassing a clear computational strategy that yields micromorphic stress-strain states from composite microstructures. The method has been used successfully in an FE$^2$ setting~\cite{biswas2019nonlinear,zhi2022} and applied to tetra-chiral structures~\cite{biswas2020micromorphic}.
		  
		\item Lastly, the method proposed by~\citet{rokovs2019} for structures involving pattern transformations deserves mention. The authors incorporate a predefined fluctuation field in the kinematic ansatz, considering the main microstructural deformation modes. As a result, a homogenized micromorphic continuum naturally emerges at the macroscale.
		
		\end{itemize}
		 
		\subsection{From discrete material systems to generalized continua}  In addition to the homogenization schemes described above, it is important to note that the homogenization of microstructured materials with inherent rotational degrees of freedom naturally yields effective micropolar quantities, e.g., in cellular and metamaterials~\cite{chen1998fracture} or granular media~\cite{ehlers2003particle}. Indeed, homogenized micropolar stress-strain states were applied recently in data-driven computing~\cite{karapiperis2021a}. Furthermore, homogenization from discrete systems has also been shown to yield micromorphic effective states~\cite{regueiro2014micromorphic,misra2015,misra2020,ehlers2020}, which can be employed for data generation in related applications.
		
		\subsection{Data-driven identification} 
		Another interesting opportunity arising from the present work, in the context of data sampling, is to extend the recent data-driven identification framework of~\citet{stainier2019model} to micromorphic continua (or subcases thereof). In this setting, a data-driven problem supplied with generalized strain data and enforced to comply with the associated balance equations is expected to yield the corresponding generalized stresses. Therefore, a framework of this type is an excellent candidate for data acquisition.

            \medskip
            
		In summary, several opportunities exist for extracting generalized stress-strain states from microstructured material responses, allowing the generation of the required datasets for our current data-driven scheme. Therefore, future developments should focus on implementing homogenization schemes or data-driven identification methods for data acquisition in applications specific to material systems.
}

\small

\end{document}